\documentclass[psamsfonts]{conm-p-l}

\usepackage{amsmath,amssymb}
\newtheorem{theorem}{Theorem}[section]
\newtheorem{lemm}[theorem]{Lemma}
\newtheorem{prop}[theorem]{Proposition}

\theoremstyle{definition}
\newtheorem{defi}[theorem]{Definition}
\newtheorem{example}[theorem]{Example}

\newtheorem{coro}[theorem]{Corollary}
\theoremstyle{remark}
\newtheorem{remark}[theorem]{Remark}

\newcommand{\ad}{\mathrm{ad}}
\allowdisplaybreaks

\numberwithin{equation}{section}

\def\lg{\langle}
\def\rg{\rangle}

\def\pa{\partial}

\def\al{\alpha}

\def\ep{\epsilon}

\def\p{\partial}

\def\z{\mathbb{Z}}

\def\dim{\hbox{dim}}
\def\ad{\hbox{ad}}

\def\a{\alpha}

\def\mod{\hbox{mod}}

\newfont{\df}{eufm10}
\def\vep{\varepsilon}
\def\ep{\epsilon}

\def\dim{\hbox{\rm dim}\,}

\def\ad{\hbox{\rm ad}\,}

\begin{document}

\title[Modular Quantizations of Lie Algebras of Cartan type $H$]
{Modular Quantizations of Lie Algebras\\ of Cartan Type $H$ via Drinfeld Twists}

\author[Tong]{Zhaojia Tong}
\address{Department of Mathematics, Shanghai University,
Baoshan Campus, Shangda Road 99, Shanghai 200444, PR China}
\email{tongzhaojia@gmail.com}

\thanks{$^\star$N.H., corresponding author, supported in part by the NNSF of China (No. 11271131)}
\thanks{$^*$X.L.,  supported by the NNSF of China (No. 10901085)}

\author[Hu]{Naihong Hu$^\star$}
\address{Department of Mathematics, Shanghai Key Laboratory of Pure Mathematics
and Mathematical Practice, East China Normal University, Minhang Campus, Dong Chuan
Road 500, Shanghai 200241, PR China}
\email{nhhu@math.ecnu.edu.cn}

\author[Wang]{Xiuling Wang$^*$}
\address{School of Mathematical Sciences and  LPMC, Nankai University,
Tianjin  300071, PR China}\email{xiulingwang@nankai.edu.cn}

\subjclass[2010]{Primary 17B37, 17B62; Secondary 17B50}

\dedicatory{Dedicated to Professor Helmut Strade in honor of his 70th
birthday}

\keywords{Lie bialgebra, $r$-matrix, Drinfel'd twist, modular quantization,
Lie algebra of generalized Cartan type $H$, restricted Hamiltonian algebra,
Hopf algebra of prime-power dimension}
\begin{abstract}
We construct explicit Drinfel'd twists for the Lie algebras of generalized
Cartan type $H$ in characteristic $0$ and also obtain the corresponding
quantizations and their integral forms. By using modular reduction and
base changes, we derive certain quantizations of the restricted universal
enveloping algebra $\mathbf u(\mathbf{H}(2n;\underline{1}))$ of the
restricted Hamiltonian algebra $\mathbf{H}(2n;\underline{1})$ in prime
characteristic $p$. These quantizations are new non-pointed Hopf algebras
of prime-power dimension $p^{p^{2n}-1}$ and contain the well-known
Radford algebras as Hopf subalgebras. As a by-product we also obtain
some Jordanian quantizations of $\mathfrak {sp}_{2n}$.
\end{abstract}

\maketitle

This paper is a continuation of \cite{HW1} and \cite{HW2} in which modular
quantizations of Lie algebras of Cartan types $W$ and $S$ were
studied. In the present paper we consider the same questions, both
for the Lie algebras of generalized Cartan type $H$ in characteristic
$0$ (for the definition, see \cite{OZ}) and for the restricted Hamiltonian
algebras $\mathbf{H}(2n;\underline{1})$ in the modular case (for
the definition, see \cite{H} and \cite{HR}).

For the convenience of the reader we review some previous related work.
In \cite{Dr1}, Drinfel'd raised the question of the existence of universal
quantizations for Lie bialgebras. In \cite{EK1} and \cite{EK2} Etingof and Kazhdan
gave a positive answer to this question for Lie bialgebras coming from finite-
and infinite-dimensional Lie algebras defined by generalized Cartan matrices.
Later, Enriquez-Halbout \cite{EH} showed that, in principle, any coboundary
Lie bialgebra can be quantized via a certain Etingof-Kazhdan quantization
functor, and Geer \cite{G} extended the work of Etingof and Kazhdan from
Lie bialgebras to the setting of Lie superbialgebras. In view of this, it is natural
to consider the quantizations of Lie algebras of Cartan type which are
defined by differential operators. In 2004 Grunspan \cite{CG} obtained the
quantization of the infinite-dimensional Witt algebra $\bold W$ in characteristic
$0$ by using the twist found by Giaquinto and Zhang in \cite{AJ}, but his approach
didn't work for the quantum version of the finite-dimensional Witt algebra $\bold
W(1;\underline{1})$ in characteristic $p>0$. The second and the third author
obtained in \cite{HW1} quantizations of the generalized Witt algebra $\bold W$
in characteristic $0$ and of the Jacobson-Witt algebras $\bold W(n;\underline{1})$
in characteristic $p>0$. These quantizations are new families of noncommutative
and noncocommutative Hopf algebras of dimension $p^{1+np^n}$ in prime characteristic
$p$, while in the rank $1$ case, \cite{HW1} recovered Grunspan's work in characteristic
$0$ and gave the desired quantum version in characteristic $p>0$.

Although, in principle, the possibility to quantize an arbitrary Lie
bialgebra has been proved (\cite{EK1}, \cite{EK2}, \cite{ES},
\cite{EH}, and \cite{G}), it seems difficult to obtain explicit formulas
for the Hopf algebra operations. In particular, only a few kinds of
twists with explicit expressions for the Hopf algebra operations
are known (see \cite{NR}, \cite{VO}, \cite{AJ}, \cite{KL}, \cite{KLS},
and the references therein). In this paper we start with an explicit
Drinfel'd twist considered in \cite{AJ} and \cite{CG}. In fact, this Drinfel'd twist  is
essentially a variation of the Jordanian twist (see the proof in \cite{HW2})
which first (but in a different form) appeared in a
paper by Coll, Gerstenhaber, and Giaquinto \cite{CGG}, and
was also employed extensively by Kulish et al (see \cite{KL}, \cite{KLS}, etc.).
By using such explicit Drinfel'd twists, we obtain vertical basic twists
and horizontal basic twists for the Lie algebras of generalized
Cartan type $H$ and the corresponding quantizations in
characteristic $0$. As in types $W$ and $S$ one can obtain from these
basic twists many other Drinfel'd twists. For the modular case
we first have to find the integral forms of the quantizations for
the Lie algebras of generalized Cartan type $H$. The crucial observation
is that we have to work over the so-called ``positive" part subalgebra
$\mathbf{H}^+$ of the Lie algebra $\mathbf{H}$ of generalized Cartan
type $H$. It is an infinite-dimensional simple Lie algebra defined over a
field of characteristic $0$, whereas over a field of characteristic $p$, it
contains a maximal ideal $J_{\underline{1}}$ and the corresponding
quotient is exactly the algebra $\mathbf{H}'(2n;\underline{1})$. Its
derived subalgebra $\mathbf{H}(2n;\underline{1})=\mathbf{H}'(2n;
\underline{1})^{(1)}$ is a restricted simple Lie algebra of Cartan type
$H$. Secondly, in order to obtain the desired finite-dimensional
quantizations of the restricted universal enveloping algebra of the
Hamiltonian algebra $\mathbf{H}(2n;\underline{1})$, we need to
carry out a modular reduction process: reduction modulo $p$ of an
integral form of the universal enveloping algebra and another reduction
to obtain deformations of the restricted universal enveloping algebra,
and accordingly suitable base changes. These are the other
two crucial technical points of this paper. As a result we obtain a new class of
noncommutative and noncocommutative Hopf algebras whose dimension is a
power of the characteristic of the ground field.

The paper is organized as follows. In Section 1, we recall some
definitions and basic facts  related to  the Lie algebras of Cartan
type $H$ and Drinfel'd twists. In Section 2,
we construct some  Drinfel'd twists  for  the Lie
algebras  of generalized Cartan type $H$, including
the vertical basic twists and the horizontal basic twists.
In Section 3, we explicitly quantize Lie bialgebra
structures of the Lie algebras of generalized Cartan type $H$ by
the vertical basic Drinfel'd twists, and by using
similar methods as in type $S$, we obtain the modular quantizations
of the restricted universal enveloping algebra of the Hamiltonian algebra
$\mathbf{H}(2n;\underline{1})$. In Section 4, we use horizontal basic
twists to get some new modular quantizations of horizontal type
of $\mathbf u(\mathbf{H}(2n;\underline{1}))$ which contain some
modular quantizations of the Lie algebra $\mathfrak{sp}_{2n}$ derived from  the
Jordanian twists (cf.~\cite{KLS}). Finally, we present some open questions.

\section{Preliminaries}
\subsection{The generalized Cartan-type Lie algebra $\textbf{H}$ and its subalgebra $\textbf{H}^+$}
We recall the definition of the generalized Cartan-type Lie
algebra $\textbf{H}$ from \cite{OZ} and some basic facts about its structure.

Let $\mathbb{F}$ be a field with char$(\mathbb{F})=0$ and let $n >
 0$. Let $\mathbb{Q}_{2n}:=\mathbb{F}[x^{\pm 1}_{-n}, \ldots,x^{\pm 1}_{-1},
x^{\pm 1}_{1}$, $\ldots, x^{\pm
 1}_{n}]$ be a Laurent polynomial algebra and let
 $\partial_i$ denote the degree operator $x_i \frac{\partial}{\partial
 x_i}$. Set $T:=\bigoplus\limits_{i=1}^{n} (\mathbb{Z}
 \partial_{-i}\oplus\mathbb{Z} \partial_i)$, and set $x^{\alpha}:=x_{-n}^{\alpha_{-n}}  \cdots x_{-1}^{\alpha_{-1}} x_1^{\alpha_1} \cdots x_n^{\alpha_n}$
 for $\alpha=(\alpha_{-n},{\ldots},\alpha_{-1}, \alpha_1, {\ldots}, \alpha_n)\in \mathbb{Z}^{2n}$. In particular,
 $x_i=x^{\ep_i}, \ep_i=(\delta_{-n,i},{\ldots},\delta_{-1,i},
 \delta_{1,i}$, ${\ldots},\delta_{n,i})$.
  We can define a bilinear map
\begin{gather*}
\langle \cdot,\cdot\rangle: T \times \mathbb{Z}^{2n}
\longrightarrow
\mathbb{Z}\\
\langle\p,\alpha \rangle \mapsto \sum\limits_{i=1}^{n}(a_{-i}\alpha_{-i}
+ a_i \alpha_i),
\end{gather*}
for $\partial=\sum\limits_{i=1}^{n} (a_{-i}\partial_{-i}+a_i \partial_i)\in T$ and $\alpha=(\alpha_{-n},{\ldots},\alpha_{-1},\alpha_1,{\ldots},\alpha_n)\in
\mathbb{Z}^{2n}$. It is easy to see that this bilinear map is
non-degenerate in the sense that
\[
\left.
\begin{array}{ccccccccc}
 \partial (\alpha)&=& \langle \partial, \alpha \rangle &=&0 \quad(\forall
~\partial &\in &T)\,&\Longrightarrow &\alpha=\underline{0}\,,\\
\partial (\alpha)&=&\langle \partial,\alpha \rangle &=&0
\quad(\forall~\alpha &\in &\mathbb{Z}^{2n})\,& \Longrightarrow &\, \partial=0\,,
\end{array}
\right.
\]
where $\underline{0}=(0,{\ldots},0)$. For later use, we set $\underline{1}=(1,{\ldots},1)$.

Define a linear map
\begin{eqnarray*}
 D_H:&~\mathbb{Q}_{2n} \longrightarrow  &\mathrm{Der}(\mathbb{Q}_{2n})~~~~~~~~~~~~~~~~~~~~~~\\
     & x^{\alpha}      \longmapsto      & \sum\limits_{i=1}^{n}
 x^{\alpha -\ep_{-i}-\ep_{i}}(\p_{-i}(\alpha) \p_i
-\partial_i(\alpha) \p_{-i})\\
&&=\sum\limits_{i=1}^{n} x^{\alpha
-\ep_{-i}-\ep_{i}}(\alpha_{-i}\p_{i} -\alpha_i \p_{-i}).
\end{eqnarray*}
We can see that the kernel of this map is $\mathbb{F}$. The image $D_H(\mathbb{Q}_{2n})$ of
this map is a Lie algebra under the bracket
\begin{eqnarray*}
&[D_H (x^{\alpha}),~D_H (x^{\beta})]&=D_H\left(\sum\limits_{i=1}^{n}\Big(\frac{\p
x^{\alpha}}{\p x_{-i}}
\frac{\p x^{\beta}}{\p x_i}-\frac{\p x^{\alpha}}{\p x_i}\frac{\p x^{\beta}}{\p x_{-i}}\Big)\right)~~~~~~~~~~\\
&&=D_H\Big(\sum\limits_{i=1}^{n} \big(\p_{-i}(\alpha)
\p_i(\beta)-\p_{i}(\alpha) \p_{-i}(\beta)\big)x^{\alpha +\beta -\ep_i
-\ep_{-i}}\Big)\\
&&=D_H\Big(\sum\limits_{i=1}^{n}\big(\alpha_{-i}\beta_i-\alpha_i
\beta_{-i}\big)x^{\alpha +\beta -\ep_i -\ep_{-i}}\Big).
\end{eqnarray*}
The derived algebra
$\mathbf{H}=[D_H(\mathbb{Q}_{2n}), D_H(\mathbb{Q}_{2n})]$ is the
Lie algebra of generalized Cartan type $H$,  which has codimension $1$ in
$D_H(\mathbb{Q}_{2n})$ and is known to be a simple algebra. Moreover, $\{D_H
(x^{\alpha}) \mid \alpha\in \mathbb{Z}^{2n} \setminus \{\underline{0},
-\underline{1}\} \}$ is a basis of $\mathbf H$ (cf.~\cite{OZ}).

Define $D_i=\frac{\p}{\p x_i}$. For $D_H
(x^{\alpha})=\sum\limits_{i=1}^{n} x^{\alpha -\ep_i
-\ep_{-i}}(\alpha_{-i} \p_i-\alpha_i \p_{-i})$ we have
$D_H(x^{\alpha})$ $=\sum\limits_{i=1}^{n} \alpha_{-i} x^{\alpha
-\ep_{-i}} D_i- \alpha_i x^{\alpha-\ep_i} D_{-i}$. Let $\mathcal{K}$
be an arbitrary field and set
$\mathbf{H}^+=\mathrm{Span}_{\mathcal{K}}\{ D_H (x^{\alpha}) \mid
\alpha ~\in \mathbb{Z}_+^{2n}\setminus
\{\underline{0}\}  \}$, which via the identification
$x^{\alpha}D_i$ with $x^{\alpha - \ep_i}\p_i$ ($\alpha-\ep_i
\in  \mathbb{Z}^{2n}$) can be considered as a Lie subalgebra (namely, the
``positive" part) of the Lie algebra of generalized Cartan type $H$
over $\mathcal{K}$.

\subsection{The Hamiltonian algebra $\mathbf{H}(2n;\underline{1})$}
Assume now that char$(\mathcal{K})=p$. Then by definition the
Jacobson-Witt algebra $\mathbf{W}(2n;\underline{1})$
is a restricted simple Lie algebra over $\mathcal{K}$ (see \cite{HR}). Its
$p$-Lie algebra structure is given by $D^{[p]}=D^p,\; \forall\, D
\in \mathbf{W}(2n;\underline{1})$ and $\{\,x^{(\alpha)}D_j
\mid -n \leq j\leq n,\, j\neq 0, \ 0 \leq \alpha \leq \tau \}$ is a basis,
where $\tau=(p{-}1,{\ldots},p{-}1)$,
$\ep_i=(\delta_{-n,i},{\ldots},\delta_{-1,i}
,\delta_{1,i},{\ldots},\delta_{n,i})$ with $x^{(\ep_i)}=x_i$,
$x^{(\al)}$ $\in\mathcal {O}(2n;\underline{1})=\text{Span}_{\mathcal K}\{\, x^{(\alpha)} \mid 0 \leq \alpha
\leq \tau \}$. The latter is the restricted divided power algebra with
multiplication $x^{(\alpha)} x^{(\beta)}=\binom{\alpha+\beta}{\alpha} x^{(\alpha
+\beta)}$,  where
 $\binom{\alpha +\beta }{\alpha}=
 \prod\limits_{i=1}^{n}\binom{\alpha_{-i}+\beta_{-i}}{\alpha_{-i}}\binom{\alpha_i + \beta_i}{\alpha_i}
$, and the convention that $x^{(\alpha)}=0$ if $\alpha$ has a
component $\alpha_j
 < 0$ or $\geq p$.
Note that
$\mathcal{O}(2n;\underline{1})_i$ $=\mathrm{Span}_\mathcal{K} \{
x^{(\alpha)} \mid 0 \leq \alpha \leq \tau, \ |\alpha|=i \}$ where
$|\alpha|=\sum\limits_{i=1}^{n} (a_{-i}{+}a_i)$. Moreover, ${\bf
W}(2n;\underline{1})=\mathrm{Der}_{\mathcal{K}}(
\mathcal{O}(2n;\underline{1}))$ and inherits a gradation from
$\mathcal{O}(2n;\underline{1})$ by means of ${\bf
W}(2n;\underline{1})_i=\sum\limits_{j=1}^{n} (
\mathcal{O}(2n; \underline{1})_{i+1} D_{-j}+\mathcal{O}(2n; \underline{1})_{i+1}
D_j)$. Define $D_H:
\mathcal{O}(2n; \underline{1})\rightarrow {\bf W}(2n;\underline{1})$
as $D_{H}(x^{(\alpha)})=\sum\limits_{i=1}^{n}
(x^{(\alpha-\epsilon_{-i})}D_i-x^{(\alpha-\epsilon_{i})}D_{-i})$.
Then the subspace $\mathbf{H}'(2n; \underline{1})$ $:=D_H(\mathcal{O}(2n; \underline{1}))$ is a subalgebra of ${\bf
W}(2n; \underline{1})$. Its derived algebra
$\mathbf{H}(2n; \underline{1})$ is called the Hamiltonian algebra,
$\mathbf{H}(2n; \underline{1})=\bigoplus\limits_{i=-1}^{s}
\mathbf{H}(2n; \underline{1}) \bigcap$  ${\bf
W}(2n;\underline{1})_i$ is graded with $s=|\tau|-3$. Then by
Proposition 4.4.4 and Theorem 4.4.5 in \cite{HR},
$\mathbf{H}(2n; \underline{1})=\mathrm{Span}_{\mathcal{K}}\{ D_H
(x^{(\alpha)}) \mid  x^{(\alpha)} \in \mathcal{O}(2n; \underline{1}),
\ 0 \leq \alpha<\tau\}$ is a $p$-subalgebra of
${\bf W}(2n;\underline{1})$ with restricted gradation.

By definition (cf. \cite{HR}), the restricted universal enveloping
algebra ${\bf u(H}(2n;\underline{1}))$ is isomorphic to $U({\bf
H}(2n;\underline{1}))/I$, where $I$ is the Hopf ideal of $U({\bf
H}(2n;\underline{1}))$ generated by $(D_H (x^{(\ep_i +\ep_{-i})}))^p
- D_H (x^{(\ep_i+\ep_{-i})})$, $(D_H (x^{(\alpha)}))^p$ with $\alpha
\neq \ep_i +\ep_{-i}$ for $ 1 \leq i \leq n$. Since
$\mathrm{dim}_\mathcal{K} {\bf H}(2n;\underline{1})=p^{2n}-2$, we
have $\mathrm{dim}_\mathcal{K}{\bf{u}}({\bf{
H}}(2n;\underline{1}))=p^{p^{2n}-2}$.

\subsection{Quantization by Drinfel'd twists}
The following result is well known (see \cite{CP}).
\begin{lemm}{ \label{twist1}
Let $(A,m,\iota,\Delta_0,\vep,S_0)$ be a Hopf algebra over a
commutative ring. A Drinfel'd twist $\mathcal {F}$ on $A$ is an
invertible element of $A \otimes A$ such that
\begin{gather*}
(\mathcal{F} \otimes 1)(\Delta_0 \otimes \text{\rm
Id})(\mathcal{F})=(1
\otimes \mathcal{F})(\text{\rm Id} \otimes \Delta_0)(\mathcal{F}),\\
(\varepsilon \otimes \text{\rm Id})(\mathcal{F})=1=(\text{\rm Id} \otimes
\varepsilon)(\mathcal{F}).
\end{gather*}

Then $w=m(\text{\rm Id} \otimes S_0)(\mathcal{F})$ is invertible
in $A$ with $w^{-1}=m(S_0 \otimes \text{\rm Id})(\mathcal{F}^{-1})$.

Moreover, if we define $\Delta: A \longrightarrow A \otimes
A$ and $S: A \longrightarrow A$ by
$$\Delta (a):=\mathcal{F} \Delta_0 \mathcal{F}^{-1},~~~~S(a):=w S_0 (a)
w^{-1},$$
then $(A,m,\iota,\Delta,\vep,S)$ is a new Hopf algebra, called
the twisting of A by the Drinfel'd twist $\mathcal{F}$.}
\end{lemm}

Let $\mathbb{F}[[t]]$ be a ring of formal power series over a
field $\mathbb{F}$ with char\,$\mathbb{F}=0$. Assume that $L$ is a
triangular Lie bialgebra over $\mathbb{F}$ with a classical
$r$-matrix $r$ (see
\cite{D} and \cite{ES}). Let $U(L)$ denote the universal enveloping
algebra of $L$ with the standard Hopf algebra structure
$(U(L),m,\iota,\Delta_0,\vep,S_0)$.

Let us consider the topologically free $\mathbb{F}[[t]]$-algebra
$U(L)[[t]]$ (for  the definition, see p.~4 of \cite{ES}), which can
be viewed as an associative $\mathbb{F}$-algebra of formal power
series with coefficients in $U(L)$. Naturally, $U(L)[[t]]$ has an Hopf
algebra structure induced from the one on $U(L)$.
By abuse of notation, we denote it by
$(U(L)[[t]],m,\iota,\Delta_0,\vep, S_0)$.

\begin{defi}{\upshape \cite{HW1}
For a triangular Lie
bialgebra $L$ over $\mathbb{F}$ with char($\mathbb{F}$)$=0$,
$U(L)[[t]]$ is called a \textit{quantization} of $U(L)$ by a
Drinfel'd
twist $\mathcal{F}$ over $U(L)[[t]]$ if
$U(L)[[t]]/tU(L)[[t]] \cong U(L)$ as algebras, and $\mathcal{F}$ is determined
by its $r$-matrix $r$ (namely, its Lie bialgebra structure).}
\end{defi}

\section{Drinfel'd twists over $U(\mathbf{H})[[t]]$}
\subsection{Construction of Drinfel'd twists}
Let $L$ be a Lie algebra containing linearly independent elements
$h$ and $e$ satisfying $[h,e]=e$. Then the classical
$r$-matrix $r=h \otimes e -e \otimes h$ equips $L$ with the
structure of a triangular coboundary Lie bialgebra (see \cite{M}).
In order to obtain an explicit description of a quantization of $U(L)$
by a Drinfel'd twist $\mathcal{F}$ over $U(L)[[t]]$, we need an
explicit construction for such a Drinfel'd twist. In what follows, we
shall see that such a Drinfel'd twist depends on the choice of two distinguished elements
$h$  and $e$ arising from its $r$-matrix $r$.

For any element of a unital $R$-algebra (where $R$ is any ring) and $a \in R$,
we set
\begin{gather*}
x_{a}^{\langle m \rangle}:=(x+a)(x+a+1)\cdots(x+a+m-1),\\
x_{a}^{[m]}:=(x+a)(x+a-1)\cdots (x+a-m+1),
\end{gather*}
and then put $x^{\langle m \rangle}:=x_{0}^{\langle m \rangle}\,,\,
x^{[m]}:=x_{0}^{[m]}$.

Note that $h$ and $e$ satisfy the following identities
\begin{eqnarray*}
e^s\cdot h_a^{[m]}=h_{a-s}^{[m]}\cdot e^s, \\
e^s\cdot h_a^{\langle m \rangle }=h_{a-s}^{\langle m \rangle }\cdot
e^s,
\end{eqnarray*} where $m$ and $s$ are non-negative
integers, and $a \in \mathbb{F}$.

For $a \in \mathbb{F}$ we set
$\mathcal{F}_a:=\sum\limits_{r=0}^{\infty}\frac{(-1)^r}{r!}h_a^{[r]}\otimes
e^rt^r, F_a:=\sum\limits_{r=0}^{\infty}\frac{1}{r!}h_a^{\lg
r\rg}\otimes e^rt^r,$ $ u_a:=m\cdot(S_0\otimes \text{\rm Id})(F_a),
v_a:=m\cdot(\text{\rm Id}\otimes S_0)(\mathcal{F}_a).$ Write
$\mathcal{F}:=\mathcal{F}_0,\, F:=F_0,\,u:=u_0,\,v:=v_0$. Since
$S_0(h_a^{\lg r\rg})=(-1)^rh_{-a}^{[r]}$ and $S_0(e^r)=(-1)^re^r$,
one has $ v_a=\sum\limits_{r=0}^{\infty}\frac{1}{r!}h_a^{[r]}
e^rt^r, \quad
u_b=\sum\limits_{r=0}^{\infty}\frac{(-1)^r}{r!}h_{-b}^{[r]} e^rt^r.
$

\begin{lemm} \label{1.6}$($\cite{CG}$)$
For $a,b \in \mathbb{F}$ one has
$$
\mathcal{F}_a F_b=1\otimes(1-et)^{a-b} \quad\text{and }\quad v_a
u_b=(1-et)^{-(a+b)}.
$$
\end{lemm}

\begin{coro}\label{1.7}{
For $a \in \mathbb{F}$, $\mathcal{F}_a$ and $u_a$ are invertible
with $\mathcal{F}_a^{-1}=F_a$ and $u_a^{-1}=v_{-a}$.  In particular,
$\mathcal{F}^{-1}=F$ and $u^{-1}=v$.}
\end{coro}

\begin{lemm}\label{1.8} {
$($\cite{HW1}$)$ For every positive integer $r$ we have
$$\Delta_0(h^{[r]})=\sum\limits_{i=0}^r \dbinom{r}{i}h^{[i]}\otimes
h^{[r-i]}.$$  Furthermore, $\Delta_0(h^{[r]})=\sum\limits_{i=0}^r
\dbinom{r}{i}h^{[i]}_{-s}\otimes h^{[r-i]}_s$ for any $s \in
\mathbb{F}$.}
\end{lemm}

\begin{prop} \label{twist}$($\cite{CG},\cite{HW1}$)$
{If a Lie algebra $L$ contains a two-dimensional solvable
Lie subalgebra with a basis $\{h, e\}$ satisfying $[h, e]=e$, then
$$\mathcal{F}=\sum\limits_{r=0}^{\infty}\frac{(-1)^r}{r!}h^{[r]}
\otimes\, e^rt^r$$ is a Drinfel'd twist on $U(L)[[t]]$.}
\end{prop}

\subsection{Basic Drinfel'd twists}
Take two distinguished elements $e:=D_H (x^{\alpha})$ and $h:=D_H
(x^{\ep_{-i}+\ep_i})$ such that $[h,e]=e$, where $ 1 \leq i \leq n$.
It is easy to see that then $\alpha_i -\alpha_{-i}=1$. The next result follows from the main
result of \cite{M}.

\begin{prop}\label{rmatrix}{
There is a triangular Lie bialgebra structure on $\mathbf{H}$
given by the classical $r$-matrix
$$
r:=D_H (x^{\ep_{-i}+\ep_i}) \otimes D_H (x^{\alpha})-D_H (x^{\alpha})
\otimes D_H (x^{\ep_{-i}+\ep_i}),\quad  \  1 \leq i
\leq n,
$$
where $\alpha \in \mathbb{Z}^{2n}$ with $\alpha_i - \alpha_{-i}
=1,$ and $[D_H (x^{\ep_{-i}+\ep_i}),D_H (x^{\alpha})]=D_H
(x^{\alpha})$.}
\end{prop}

Fix two distinguished elements $h:=D_H (x^{\ep_{-i}+\ep_i})$ and $e:=D_H
(x^{\alpha})$ with $\alpha_i-\alpha_{-i}=1$. Then
$\mathcal{F}=\sum\limits_{r=0}^{\infty}\frac{(-1)^r}{r!}h^{[r]}\otimes
e^rt^r$ is a Drinfel'd twist on $U(\mathbf{H})[[t]]$. But the
coefficients of the quantizations of the standard Hopf algebra structure
$(U(\mathbf{H})[[t]],m,\iota,\Delta_0$, $S_0,\varepsilon)$ by
$\mathcal{F}$ may be not integral. In order to get integral forms
of these quantizations, one needs to find sufficient and necessary conditions
for the coefficients of $\mathcal{F}$ being integers.

\begin{lemm} $($\cite{CG}$)${
For any $a,\,k,\,\ell\in\mathbb{Z}$, the rational number
$a^\ell\prod\limits_{j=0}^{\ell-1}(k{+}ja)/\ell!$ is an
integer.\hfill\qed }
\end{lemm}

In view of this lemma, we are interested in the following two simple
cases:

\smallskip
$\mathrm{(i)} \ h=D_H (x^{\ep_{-k}+\ep_k}), ~e=D_H (x^{\ep_{-k}+2
\ep_{k}}), \quad (1 \leq k \leq n)$;

$\mathrm{(ii)} \ h=D_H (x^{\ep_{-k}+\ep_k}), ~e=D_H (x^{\ep_k +
\ep_m}), \quad (m \neq k, -k)$.

\smallskip
Let $\mathcal{F}(k)$ denote the corresponding Drinfel'd twist in case (i)
and $\mathcal{F}(k;m)$ denote the corresponding Drinfel'd twist in
case (ii).

\begin{defi}\label{VH twist} The twist $\mathcal{F}(k)$ $(1\leq k \leq
n)$ is called a {\it vertical basic Drinfel'd twist} and the twist $\mathcal{F}(k;m)$
$(1\leq k, m \leq n, m \neq k, -k)$ is called a {\it horizontal basic Drinfel'd
twist}.
\end{defi}

\begin{remark}\label{2.3}
In case (i) we get $n$ vertical basic Drinfel'd twists
$\mathcal{F}(1)$, $\ldots,\mathcal{F}(n)$ for
 $U(\mathbf{H}_{\mathbb{Z}}^+)[[t]]$. It is interesting to consider
the products of some of these vertical basic
Drinfel'd twists. In this way one can obtain more Drinfel'd twists which
will lead to new quantizations not only over
$U(\mathbf{H}_{\mathbb{Z}}^+)[[t]]$, but via our modular
reduction approach developed in Section 3, also over
$\mathbf{u}_{t,q}(\mathbf{H}(2n;\underline{1}))$.

In case (ii) we get $2n(n-1)$ horizontal basic Drinfel'd twists
$\mathcal F(k;m)$ over $U(\mathbf{H}^+_{\mathbb{Z}})[[t]]$. We will discuss these
twists and the corresponding quantizations in Section 4.
\end{remark}

\subsection{More Drinfel'd twists}
We consider the products of pairwise different and mutually
commutative vertical basic Drinfel'd twists and from this we can get new
quantizations not only over
$U(\mathbf{H}_{\mathbb{Z}}^+)[[t]]$ but over
$\mathbf{u}(\mathbf{H}(2n;\underline{1}))$ as well. Note that
$[\mathcal{F}(k), \mathcal{F}(k^{\prime})]=0$ for $ 1 \leq k \neq
k^{\prime} \leq n$. According to the definition of
$\mathcal{F}(k)$, this fact implies the following commutativity relations for
$ 1 \leq k \neq k^{\prime} \leq n$:
\begin{equation}\label{relation0}
\begin{split}
(\mathcal{F}(k)\otimes 1)(\Delta_0\otimes\text{\rm Id})
(\mathcal{F}(k^{\prime}))&=(\Delta_0\otimes\text{\rm Id})
(\mathcal{F}(k))(\mathcal{F}(k^{\prime})\otimes 1),\\
(1\otimes \mathcal{F}(k))(\text{\rm Id}\otimes\Delta_0)
(\mathcal{F}(k^{\prime}))&=(\text{\rm Id}\otimes\Delta_0)
(\mathcal{F}(k))(1\otimes\mathcal{F}(k^{\prime})),
\end{split}
\end{equation}
which give rise to the following property:

\begin{theorem}\label{twist3}
$\mathcal{F}(k)\mathcal{F}(k^\prime)$ $(1 \leq k \neq k' \leq n )$ is
also a Drinfel'd twist on $U(\mathbf{H}_{\mathbb{Z}}^+)[[t]]$.
\end{theorem}
\begin{proof}
Note that $\Delta_0\otimes\text{\rm id}$, $\text{\rm
id}\otimes\Delta_0$, $\varepsilon_0\otimes\text{\rm id}$ and
$\text{\rm id}\otimes\varepsilon_0$ are algebraic homomorphisms.
According to Lemma \ref{twist1}, it suffices to check that
$$\big(\mathcal{F}(k)\mathcal{F}(k^\prime)\otimes
1\big)(\Delta_0\otimes
\text{\rm Id})\big(\mathcal{F}(k)\mathcal{F}(k^\prime)\big)\\
= \big(1\otimes \mathcal{F}(k)\mathcal{F}(k^\prime)\big)(\text{\rm
Id}\otimes\Delta_0) \big(\mathcal{F}(k)\mathcal{F}(k^\prime)\big).$$

By using $(\ref{relation0})$, we have
\begin{equation*}
\begin{split}
\text{LHS}&=(\mathcal{F}(k)\otimes 1)(\mathcal{F}(k^\prime)\otimes
1)(\Delta_0\otimes \text{\rm Id})(\mathcal{F}(k))(\Delta_0\otimes
\text{\rm
Id})(\mathcal{F}(k^\prime))\\
&=(\mathcal{F}(k)\otimes 1)(\Delta_0\otimes \text{\rm
Id})(\mathcal{F}(k))(\mathcal{F}(k^\prime) \otimes
1)(\Delta_0\otimes \text{\rm
Id})(\mathcal{F}(k^\prime))\\
&=(1\otimes\mathcal{F}(k) )(\text{\rm Id}\otimes\Delta_0
)(\mathcal{F}(k))(1\otimes \mathcal{F}(k^\prime)(\text{\rm
Id}\otimes\Delta_0)(\mathcal{F}(k^\prime))\\
&=(1\otimes\mathcal{F}(k) )(1\otimes \mathcal{F}(k^\prime)(\text{\rm
Id}\otimes\Delta_0 )(\mathcal{F}(k))(\text{\rm
Id}\otimes\Delta_0)(\mathcal{F}(k^\prime))=\text{RHS}.
\end{split}
\end{equation*}
This completes the proof.
\end{proof}

More generally, we have the following
 \begin{coro}\label{2.5}
Let $\mathcal{F}(k_1),\ldots, \mathcal{F}(k_m)$ be $m$ pairwise
different vertical basic Drinfel'd twists and
 $[\mathcal{F}(k_i), \mathcal{F}(k_s)]=0$ for all $1\leq i\neq s\leq
 n$. Then $\mathcal{F}(k_1) \cdots \mathcal{F}(k_m)$
 is still a Drinfel'd twist.
\end{coro}
We set $\mathcal{F}_m:=\mathcal{F}(k_1)\cdots \mathcal{F}(k_m)$
and call $m$ its \textit{length}. We shall show that the twisted structures
given by Drinfel'd twists with different lengths are nonisomorphic.

\begin{defi}\label{com}(\cite{MG},\cite{HW2}) {\upshape
A Drinfel'd twist $\mathcal{F} \in A\otimes A$ on any Hopf algebra
$A$ is called \textit{compatible} if $\mathcal{F}$ commutes with the
coproduct $\Delta_0$.}
\end{defi}

In other words, twisting a Hopf algebra $A$ with a
compatible twist $\mathcal{F}$ gives exactly the same Hopf
algebra structure, that is, $\Delta_{\mathcal{F}}=\Delta_0$. The set of
compatible twists on $A$ thus forms a group.

\begin{lemm}\label{2.7}$($\cite{MG}$)$ {\upshape
Let $\mathcal{F} \in A\otimes A$ be a Drinfel'd twist on a Hopf
algebra $A$. Then the twisted structure induced by $\mathcal{F}$
coincides with the structure on $A$ if and only if $\mathcal{F}$ is
a compatible twist.}
\end{lemm}

\begin{lemm}\label{twist4}$(${\upshape \cite{HW2}$)$
Let $\mathcal{F}, \mathcal{G} \in A\otimes A$ be Drinfel'd twists on
a Hopf algebra $A$ with
$\mathcal{F}\mathcal{G}=\mathcal{G}\mathcal{F}$ and $\mathcal{F}\neq
\mathcal{G}$. Then $\mathcal{F}\mathcal{G}$ is a Drinfel'd twist.
 Furthermore, $\mathcal{G}$ is a Drinfel'd twist on
$A_{\mathcal{F}}$, $\mathcal{F}$ is a Drinfel'd twist on
$A_{\mathcal{G}}$
 and
$\Delta_{\mathcal{F}\mathcal{G}}=(\Delta_{\mathcal{F}})_{\mathcal{G}}
=(\Delta_{\mathcal{G}})_{\mathcal{F}}$.}
\end{lemm}

\begin{prop}\label{2.9}{\upshape
Drinfel'd twists
$\mathcal{F}^{\zeta(i)}:=\mathcal{F}(1)^{\zeta_1}\cdots\mathcal{F}(n)^{\zeta_n}$,
where
$\zeta(i)=(\zeta_1,\ldots,\zeta_{n})=(\underbrace{1,\ldots,1}_i,0,{\ldots},0)\in
\mathbb Z_2^{n}$, lead to $n$ different twisted Hopf algebra
structures on $U(\mathbf{H}_{\mathbb{Z}}^{+})[[t]]$.}
\end{prop}
\begin{proof} The proof is the same as that of Proposition 2.15 in \cite{HW2}.
\end{proof}

\section{Quantizations of vertical type for Lie bialgebras of Cartan type $H$}

In this section we explicitly quantize the Lie bialgebras of type
$\mathbf{H}$ by the vertical basic Drinfel'd twists  and obtain
certain quantizations of the restricted universal enveloping algebra
$\mathbf{u}(\mathbf{H}(2n;\underline{1}))$ by modular reduction
and base change.

\subsection{Integral quantizations of the $\mathbb{Z}$-form $\mathbf{H}_{\mathbb{Z}}^{+}$ in characteristic 0}
For the universal enveloping algebra $U(\mathbf{H})$ of the Lie
algebra $\mathbf{H}$ over $\mathbb{F}$ we denote by $(U(\mathbf{H}),m,\iota,\Delta_0,S_0,
\varepsilon)$ the standard
Hopf algebra structure. We can perform the process of twisting the
standard Hopf structure by the vertical Drinfel'd twist $\mathcal{F}(k)$ with $h:=D_H (x^{\ep_{-k}+\ep_k}),
e:=D_H(x^{\ep_{-k}+2\ep_k})$. In order to simplify the formulas, let us introduce
the operator $d^{(\ell)}$ on $U(\mathbf{H})$ defined by
$d^{(\ell)}:=\frac{1}{\ell!} (\ad e)^\ell$. By using induction on $\ell$,
we get $d^{(\ell)}(D_H (x^{\alpha}))=A_\ell D_H (x^{\alpha +
\ell \ep_k})$, where $A_\ell:=\frac{1}{\ell
!}\prod\limits_{j=0}^{\ell-1}(\alpha_k -2\alpha_{-k}+j)\in\mathbb Z$ (the latter follows from Lemma 2.6),  and set
$A_0:=1, A_{-1}:=0$.

Recall the vertical basic Drinfel'd twist of the Lie algebra of Cartan type $S$
given by $h\,{:=}\,\pa_k {-} \pa_{-k}$, $e\,{:=}\,x^{\ep_k}(\pa_k
{-}2 \pa_{-k})$ in \cite{HW2} and the vertical basic twist of the
special algebra in characteristic $p$ given by
$h\,{:=}\,D_{k,-k}(x^{(\ep_{-k}{+}\ep_k)})$,
$e\,{:=}\,2D_{k,-k}(x^{(\ep_{-k}{+}2\ep_k)})$. Note that
$\mathbf{H}$ is a Lie subalgebra of $\mathbf{S}$, and $D_{i,-i}(x^{\al})\in\mathbf{S}$ for any
$x^{\al} \in\mathbb Q_{2n}$. By virtue of the quantizations of the
Lie algebra of Cartan type $S$ in \cite{HW2} and the formulas $D_H
(x^{\alpha})=\sum\limits_{i=1}^nD_{i,-i}(x^{\al})=\sum\limits_{i=1}^{n}
x^{\alpha {-}\ep_i{-}\ep_{-i}}(\alpha_{-i} \p_i{-}\alpha_i \p_{-i})\in
\mathbf{H}\subset \mathbf{S}$, we have the following result which gives the quantization
of $U({\bf H})$ by the Drinfel'd twist $\mathcal{F}(k)$.

\begin{lemm}\label{vrelation1}
For $D_H (x^\alpha) \in U(\mathbf{H})$ the following identities hold:
\begin{gather*}
D_H (x^{\alpha}) \cdot e^m=\sum\limits_{\ell=0}^{m} (-1)^\ell
\binom{m}{\ell} e^{m-\ell} a_\ell D_H (x^{\alpha+\ell \ep_k}),\tag{\text{\rm i}}\\
\big(\ad D_H (x^\alpha)\big)^m
(e)=\prod_{j=0}^{m-1}\big((j-1)\alpha_k-(j-2)\alpha_{-k}\big)D_H
(x^{m(\alpha-\ep_k-\ep_{-k})+\ep_{-k}+2\ep_k}),\tag{\text{\rm ii}}
\end{gather*}
where $a_\ell:=\ell! A_\ell$ and $A_\ell$ is defined as before.
\end{lemm}
\begin{proof}
It is easy to get the first formula by a direct calculation. For the second one use
induction on $m$. This is true for $m=1$, since
$\ad D_H (x^\alpha) \cdot e=[D_H (x^\alpha),~D_H
(x^{2\ep_k+\ep_{-k}})]=(2\alpha_{-k}-\alpha_k)D_H (x^{\alpha+\ep_k}).$
For $m \geq 1$ we have
{\setlength{\arraycolsep}{0pt}
\begin{eqnarray*}
&&(\ad D_H (x^\alpha))^{m+1} (e) =(\ad D_H (x^\alpha))
\prod\limits_{j=0}^{m-1}\big((j-1)\alpha_{k}-(j-2)\alpha_{-k}\big)D_H
(x^{m(\alpha-\ep_k-\ep_{-k})+\ep_{-k}+2\ep_k})\\
&&=\prod\limits_{j=0}^{m-1}\big((j-1)\alpha_{k}-(j-2)\alpha_{-k}\big)\\
&&\quad \sum\limits_{i=1}^{n}\big(\alpha_{-i}(m(\alpha-\ep_k-\ep_{-k})+\ep_{-k}+2\ep_{k})_i-\alpha_i
\big(m(\alpha-\ep_{-k}-\ep_k)+\ep_{-k}+2\ep_k\big)_{-i}\big) \\
&&\qquad D_H (x^{m(\alpha-\ep_k-\ep_{-k})+\ep_{-k}+2\ep_k+\alpha-\ep_i-\ep_{-i}})\\
&&=\prod\limits_{j=0}^{m-1}\big((j-1)\alpha_{k}-(j-2)\alpha_{-k}\big)
\big( (m-1)\alpha_k-(m-2)\alpha_{-k}\big)D_H
(x^{m(\alpha-\ep_k-\ep_{-k})+\alpha+\ep_k})\\
&&=\prod\limits_{j=0}^{m}\big((j-1)\alpha_{k}-(j-2)\alpha_{-k}\big)
D_H (x^{(m+1)(\alpha-\ep_k-\ep_{-k})+2\ep_k+\ep_{-k}}).
\end{eqnarray*}}
\end{proof}
By a direct calculation and similar arguments as in Lemma 3.3 of \cite{HW2} we obtain from Lemma \ref{vrelation1}:
\begin{lemm}\label{vrelation2}
For $D_H (x^\alpha) \in U(\mathbf{H})$,  $a\in \mathbb{F}$, and $s \in \mathbb{Z}$ one has
\begin{gather*}
\big((D_H (x^\alpha))^s \otimes 1 \big) \cdot F_a=F_{a+s(\alpha_{-k}
-\alpha_{k})} \cdot \big((D_H
(x^\alpha))^s \otimes 1\big),\tag{\text{\rm i}}\\
(D_H (x^{\alpha}))^s \cdot u_a= u_{a+s(\alpha_k -\alpha_{-k})}
\Big(\sum\limits_{\ell=0}^{\infty}d^{(\ell)} (D_H (x^{\alpha}))^s
h_{1-a}^{\lg \ell \rg} t^\ell \Big),\tag{\text{\rm ii}}\\
\big (1 \otimes (D_H (x^{\alpha}))^s \big) \cdot
F_a=\sum\limits_{\ell=0}^{\infty}(-1)^\ell F_{a+\ell} \cdot
\big(h_{a}^{\lg \ell \rg} \otimes d^{(\ell)}(D_H (x^\alpha))^s t^\ell
\big).\tag{\text{\rm iii}}
\end{gather*}
\end{lemm}

\begin{theorem}\label{Hopf}
For the two distinguished elements $h:=D_H(x^{\ep_{-k}+\ep_k})$ and $e:=D_H
(x^{\ep_{-k}+2\ep_k})$ with $[h,\,e]=e$
 in the Lie algebra $\mathbf{H}$ of generalized Cartan type $H$ over
$\mathbb{F}$ there exists the structure of a noncommutative and noncocommutative Hopf algebra
$(U(\mathbf{H})[[t]],m,\iota,\Delta,S,\varepsilon)$
that leaves the algebra structure of $U(\mathbf{H})[[t]]$ undeformed and has
the following coproduct, antipode, and counit, respectively:
\begin{gather*}
\Delta(D_H (x^{\alpha}))\,{=}\,D_H (x^{\alpha}){\otimes}
(1{-}et)^{\al_k-\al_{-k}}+\sum\limits_{\ell=0}^{\infty}(-1)^\ell
h^{\lg \ell\rg}{\otimes} (1{-}et)^{-\ell}\cdot d^{(\ell)}(D_H
(x^{\al}))t^\ell,\\
S(D_H (x^{\alpha}))=-(1{-}et)^{-(\al_k
-\al_{-k})}\cdot\Bigl(\sum\limits_{\ell=0}^{\infty} d^{(\ell)}(D_H
(x^{\al}))\cdot h_1^{\lg \ell\rg}t^\ell\Bigr),
\end{gather*}
and $\varepsilon(D_H (x^{\alpha}))=0$ for any $D_H (x^{\alpha}) \in {\bf H}$.
\end{theorem}
\begin{proof}By Lemma \ref{twist1}, Lemma \ref{1.6}, and Lemma \ref{vrelation2}, we have
 \begin{eqnarray*}
\begin{split}
\Delta(D_H (x^\alpha))
&=\mathcal{F}\cdot\Delta_0(D_H (x^{\alpha}))\cdot\mathcal{F}^{-1}\\
&=\mathcal{F}\cdot(D_H (x^{\alpha}) \otimes 1)\cdot F+
\mathcal{F}\cdot(1\otimes D_H (x^{\alpha}))\cdot F
\\
&=\Bigl(\mathcal{F} F_{\al_{-k}-\al_k}\Bigr)\cdot(D_H (x^{\alpha})
{\otimes} 1)+ \sum\limits_{\ell=0}^{\infty}(-1)^\ell
\Bigl(\mathcal{F}F_{\ell}\Bigr)\cdot\Bigl(h^{\lg \ell\rg}{\otimes}
d^{(\ell)}(D_H (x^\al))t^\ell\Bigr)\\
 &=D_H (x^{\alpha}) \otimes
(1{-}et)^{\al_k-\al_{-k}}+\sum\limits_{\ell=0}^{\infty}(-1)^\ell
h^{\lg \ell\rg}\otimes (1{-}et)^{-\ell} d^{(\ell)}(D_H
(x^{\al}))t^\ell,
\end{split}
\end{eqnarray*}
\begin{eqnarray*}
\begin{split}
S(D_H (x^{\alpha}))&=u^{-1}S_0(D_H (x^{\alpha}))\,u=-v\cdot D_H
x^{\alpha} \cdot u \quad \qquad \qquad \qquad \qquad \qquad \qquad \qquad \\
 &=-v \cdot u_{\al_k-\al_{-k}}\cdot
\Bigl(\sum\limits_{\ell=0}^{\infty}d^{(\ell)}(D_H (x^\al))
h_{1}^{\lg \ell\rg}t^\ell\Bigr)
\\ &=-(1{-}et)^{(\al_{-k}-\al_{k})}\cdot\Bigl(\sum\limits_{\ell=0}^{\infty}
d^{(\ell)}(D_H (x^\al)) h_1^{\lg \ell\rg} t^\ell\Bigr).
\end{split}
\end{eqnarray*}
This completes the proof.
\end{proof}

Note that $\{ D_H (x^\alpha) \mid \alpha \in {\z}_+^{2n} \setminus
\{\underline{0}\} \}$ is a $ \z $-basis of $\mathbf{H}_{\z}^+$ and that $\mathbf{H}_{\z}^+$ is a subalgebra
of $\mathbf{H}_{\z}$ and ${\bf W}_{\z}^+$.
The following integral quantization of $\mathbf{H}_{\z}^+$ is of importance for working out quantizations of
$U(\mathbf{H}(2n;\underline{1}))$.

\begin{coro}\label{zHopf}
With the same distinguished elements $(h,\,e)$ as above, the coalgebra structure and the antipode of the
integral quantization of $U(\mathbf{H}^+_{\mathbb{Z}})$ over
$U(\mathbf{H}^+_{\mathbb{Z}})[[t]]$ by the Drinfel'd twist
$\mathcal{F}(k)$ with an undeformed algebra structure are given by
\begin{gather*}
\Delta(D_H (x^{\alpha})){=}D_H (x^{\alpha}){\otimes}
(1{-}et)^{\al_k-\al_{-k}}{+}\sum\limits_{\ell=0}^{\infty}(-1)^\ell A_\ell
h^{\lg \ell\rg}{\otimes} (1{-}et)^{-\ell}{\cdot} (D_H
(x^{\al+\ell \ep_k}))t^\ell,\\
S(D_H (x^{\alpha}))=-(1{-}et)^{-(\al_k
-\al_{-k})}\cdot\Bigl(\sum\limits_{\ell=0}^{\infty} A_\ell (D_H
(x^{\al+\ell\ep_k}))\cdot h_1^{\lg \ell\rg}t^\ell\Bigr),
\end{gather*}
and $\varepsilon(D_H (x^{\alpha}))=0$ for any $\alpha \in {\z}_+^{2n} \setminus
\{\underline{0}\}$, where $A_\ell:=\frac{1}{\ell !}\prod\limits_{j=0}^{\ell-1}(\al_k {-}2
\al_{-k}{+}j)\in\mathbb Z$ with $A_0:=1, A_{-1}:=0$ as before.
\end{coro}
\begin{proof}
  We can get the result from Theorem \ref{Hopf} and a direct calculation.
\end{proof}
\subsection{Quantizations of the Hamiltonian algebra $\mathbf{H}(2n; \underline{1})$}
In this subsection we proceed in two steps to obtain quantizations of the restricted simple
Hamiltonian algebra $\mathbf{H}(2n; \underline{1})$ in prime characteristic $p$. Firstly,
in order to obtain quantizations of the universal enveloping algebra $U(\mathbf{H}(2n;
\underline{1}))$ of $\mathbf{H}(2n;\underline{1})$, we reduce the quantization of $U
(\mathbf{H}^+_{\mathbb{Z}})$ in characteristic $0$ (see Corollary~\ref{zHopf}) modulo
$p$, and then we make a base change from $\mathcal{K}[[t]]$ to $\mathcal{K}[t]$. Secondly,
we shall apply another base change from $\mathcal{K}[t]$ to $\mathcal{K}[t]^{(q)}_p$ to the
quantization of $U(\mathbf{H}(2n;\underline{1}))$ from the first step and then a reduction
modulo an appropriate ideal to finally obtain the desired quantization of the restricted universal
enveloping algebra $\mathbf{u}(\mathbf{H}(2n;\underline{1}))$ of $\mathbf{H}(2n;\underline{1})$.

Let $\mathbb{Z}_p$ be the prime subfield of a field $\mathcal{K}$ of prime characteristic
$p$. When considering
$\mathbf{W}_{\mathbb{Z}}^+$ as a $\mathbb{Z}_p$-Lie algebra, denoted by
$\mathbf{W}_{\mathbb{Z}_p}^+$, namely,
reducing the defining relations of $\mathbf{W}_{\mathbb{Z}}^+$ modulo $p$, we see that
$(J_{\underline{1}})_{\mathbb{Z}_p}=\text{Span}_{\mathbb{Z}_p}\{x^\alpha
D_i \mid \exists\, j: \alpha_j\ge p\,\}$ is a maximal ideal of
$\mathbf{W}^+_{\mathbb{Z}_p}$ and that
$$\mathbf{W}^+_{\mathbb{Z}_p}/(J_{\underline{1}})_{\mathbb{Z}_p}
\cong \mathbf{W}(2n;\underline{1})_{\mathbb{Z}_p}
=\text{Span}_{\mathbb{Z}_p}\{x^{(\alpha)}D_{i}\mid 0\le \alpha\le
\tau, 1\le i\le 2n\}.$$ For the subalgebra
$\mathbf{H}_{\mathbb{Z}}^+$ we have
$\mathbf{H}^+_{\mathbb{Z}_p}/(\mathbf{H}^+_{\mathbb{Z}_p}\cap(J_{\underline{1}})_{\mathbb{Z}_p})
\cong \mathbf{H}'(2n;\underline{1})_{\mathbb{Z}_p}$. In the following we will denote
$\mathbf{H}^+_{\mathbb{Z}_p}\cap(J_{\underline{1}})_{\mathbb{Z}_p}$ simply
as  $(J^+_{\underline{1}})_{\mathbb{Z}_p}$.
Moreover, we have $$\mathbf{H}'(2n;\underline{1})
=\mathcal{K}\otimes_{\mathbb{Z}_p}\mathbf{H}'(2n;\underline{1})_{\mathbb{Z}_p}
=\mathcal{K}\mathbf{H}'(2n;\underline{1})_{\mathbb{Z}_p},$$ and
$\mathbf{H}^+_{\mathcal{K}}=\mathcal{K}\mathbf{H}^+_{\mathbb{Z}_p}$.

Observe that the ideal
$J^+_{\underline{1}}:=\mathcal{K}(J^+_{\underline{1}})_{\mathbb{Z}_p}$
generates an ideal of $U(\mathbf{H}^+_{\mathcal{K}})$ over
$\mathcal{K}$, denoted by
$J:=J^+_{\underline{1}}U(\mathbf{H}^+_{\mathcal{K}})$, where
$\mathbf{H}^+_{\mathcal{K}}/J^+_{\underline{1}}\cong
\mathbf{H}'(2n;\underline{1})$. It follows from the formulas of Corollary
\ref{zHopf} that $J$ is a Hopf ideal of $U(\mathbf{H}^+_{\mathcal{K}})$
satisfying $U(\mathbf{H}^+_{\mathcal{K}})/J\cong
U(\mathbf{H}'(2n;\underline{1}))$. Note that the elements $\sum a_{i,
\alpha}\frac{1}{\alpha!}D_H (x^{\alpha})$ in
$\mathbf{H}^+_{\mathcal{K}}$, where $\underline{0}\le\alpha\le\tau$, will be
identified with $\sum a_{i, \alpha}D_H (x^{(\alpha)})$ in
$\mathbf{H}'(2n;\underline{1})$ and those in $J_{\underline{1}}$
with $0$. Hence we obtain from Corollary \ref{zHopf} the quantization of
$U\big(\mathbf{H}'(2n;\underline{1})\big)$ over
$U_t(\mathbf{H}'(2n;\underline{1})):=U(\mathbf{H}'(2n;\underline{1}))\otimes_{\mathcal
K}\mathcal K[t]$ (but not necessarily over
$U(\mathbf{H}'(2n;\underline{1}))[[t]]$, see formulas (\ref{3.1}) and
(\ref{3.2}) below):

\begin{theorem}\label{Hopfprime}
For the two distinguished elements $h:=D_H (x^{(\ep_{-k}
+\ep_k)})$ and $e:=2D_H (x^{(\ep_{-k}+2\ep_k)})$ $(1 \leq k \leq n)$
the coalgebra structure and the antipode of the corresponding quantization of
$U(\mathbf{H}^\prime(2n;\underline{1}))$ over $U_t(\mathbf{H}^\prime
(2n;\underline{1}))$ with an undeformed algebra structure are given by
\begin{gather}
\Delta\big(D_H (x^{(\alpha)})\big){=}D_H (x^{(\alpha)}){\otimes}
(1{-}et)^{\alpha_k-\alpha_{-k}}{+}\sum\limits_{\ell=0}^{p-1} (-1)^\ell\overline{A}_{\ell}
h^{\lg \ell \rg}{\otimes}(1{-}et)^{-\ell} D_H
(x^{(\alpha +
\ell \ep_k)}) h_1^{\lg \ell \rg } t^\ell, \label{3.1}\\
S\big(D_H (x^{(\alpha)})\big)=-(1{-}et)^{\alpha_{-k} -\alpha_{k}}\sum\limits_{\ell=0}^{p-1} \overline{A}_\ell D_H (x^{(\alpha + \ell \ep_k)}) h_1^{\lg \ell \rg } t^\ell, \label{3.2}
\end{gather}
and $\vep\big(D_H (x^{(\alpha)})\big)=0$, where $\overline{A}_\ell:=\ell !\binom{\alpha_k {+}\ell }{\alpha_k} A_\ell
\,(\mathrm{mod}\, p)$, $A_\ell:=\frac{1}{\ell!}\prod_{j=0}^{\ell-1}
(\alpha_k {-}2\alpha_{-k}{+}j)$, $ A_0:=1$, and $A_{-1}:=0$.
\end{theorem}
\begin{proof}By Corollary \ref{zHopf}, we have
\begin{gather*}
\Delta\big(D_H (x^{(\alpha)})\big)=\frac{1}{\al !}\Delta\big(D_H (x^{\alpha})\big)\qquad \qquad\qquad\qquad\qquad\qquad\qquad\qquad\qquad\qquad\qquad\qquad\\
=D_H (x^{(\alpha)}){\otimes}
(1{-}et)^{\al_k-\al_{-k}}{+}\sum\limits_{\ell=0}^{\infty}(-1)^\ell\frac{(\al+\ell \ep_k)!}{\al !}A_\ell
h^{\lg \ell\rg}{\otimes} (1{-}et)^{-\ell}{\cdot} (D_H
(x^{(\al+\ell \ep_k)}))t^\ell\\
{=}D_H (x^{(\alpha)}){\otimes}
(1{-}et)^{\alpha_k-\alpha_{-k}}{+}\sum\limits_{\ell=0}^{p-1} (-1)^\ell\overline{A}_{\ell}
h^{\lg \ell \rg}{\otimes}(1{-}et)^{-\ell} D_H
(x^{(\alpha +
\ell \ep_k)}) h_1^{\lg \ell \rg } t^\ell.\qquad\qquad\qquad
\end{gather*}
Then we can get the other formulas by similar arguments.
\end{proof}
Note that when $\alpha{+}\ell
\ep_k=\tau$, we have $\alpha_{k}{-}2\alpha_{-k}{+}\ell{-}1=0 \,(\mod\, p)
$, i.e., $\overline{A}_\ell=0$. Thus, Theorem \ref{Hopfprime} gives a
quantization of
$U(\mathbf{H}(2n;\underline{1}))$ over $U_t(\mathbf{H}(2n;\underline{1})):=U(\mathbf{H}(2n;\underline{1}))\otimes_{\mathcal{K}}
\mathcal{K}[t]$ (which is contained in  $U(\mathbf{H}(2n;\underline{1}))[[t]]$ as a sub-Hopf algebra).
Recall that after reducing the integral quantization of $U(\mathbf{H}_{\mathbb{Z}}^{+})$ modulo $p$
we used a base change from $\mathcal{K}[[t]]$ to $\mathcal{K}[t]$. Correspondingly, the quantized object involved has been changed from
$U(\mathbf{H}(2n;\underline{1}))[[t]]$ to
$U_t(\mathbf{H}(2n;\underline{1}))$.

Denote by $I$ the ideal
of $U(\mathbf{H})(2n;\underline{1})$ over $\mathcal{K}$ generated by
$\big(D_H (x^{(\alpha)})\big)^p$ with $\alpha \neq \ep_{-i}+\ep_i$ for $\underline{0} \le \alpha <
\tau$ and
$\big( D_H (x^{(\ep_{-i}+\ep_i}))\big)^p-D_H (x^{(\ep_{-i}
+\ep_i}))$ for $1 \leq i \leq n$. Note that $\mathbf{u}\big(\mathbf{H}(2n;\underline{1})\big)=U\big(\mathbf{H}(2n;\underline{1})\big)/I$ is
of prime-power dimension $p^{p^{2n}-2}$. In order to get a
reasonable quantization of prime-power dimension
for $\mathbf{u}\big(\mathbf{H}(2n;\underline{1})\big)$ in
characteristic $p$, it is necessary to choose an appropriate underlying vector space in which the desired
$t$-deformed object would exist.
According to our modular reduction
approach, it should be induced from the
$\mathcal{K}[t]$-algebra $U_t(\mathbf{H}(2n;\underline{1}))$ in
Theorem \ref{Hopfprime}.

Firstly, we observe the following facts (for the proof, see \cite{HW2} or \cite{HW1}):

\begin{lemm}\label{relation} $(\text{\rm i})$ \ $(1-et)^p\equiv 1 \quad (\text{\rm mod}\,p, I)$.

\smallskip $(\text{\rm ii})$ \ $(1-et)^{-1}\equiv
1+et+\cdots+e^{p-1}t^{p-1} \quad (\text{\rm mod}\,
 p, I)$.

\smallskip $(\text{\rm iii})$ \ $h_a^{\lg \ell\rg} \equiv 0 \quad
(\text{\rm mod} \, p, I)$ for $\ell \geq p$, and $a\in\mathbb{Z}_p$.
\end{lemm}

Lemma \ref{relation} in conjunction with Theorem \ref{Hopfprime} indicates
that the desired $t$-defor\-mation of
$\mathbf{u}(\mathbf{H}(2n;\underline{1}))$ (if it exists) can
only occur in a certain $p$-truncated polynomial ring in one indeterminate $t$
of degree less than $p$ with coefficients in
$\mathbf{u}(\mathbf{H}(2n;\underline{1}))$, that is,
$\mathbf{u}_{t,q}(\mathbf{H}(2n;\underline{1})):=
\mathbf{u}(\mathbf{H}(2n;\underline{1}))\otimes_{\mathcal K}
\mathcal{K}[t]_p^{(q)}$ (rather than in
$\mathbf{u}(\mathbf{H}(2n;\underline{1}))$
$\otimes_{\mathcal K}\mathcal K[t]$), where $\mathcal{K}[t]_p^{(q)}$
denotes the following quotient of $\mathcal K[t]$ depending on an element
$q\in\mathcal{K}$:
$$
\mathcal{K}[t]_p^{(q)}= \mathcal{K}[t]/(t^p-qt).
\label{truncated}
$$
In order to ensure ``restrictedness" in our second reduction, it is
necessary for us to first change from $U_t(\mathbf
H(2n;\underline1))$ to $U_{t,q}(\mathbf H(2n;\underline1))$,
and then to $\mathbf u_{t,q}(\mathbf H(2n;\underline1))$ (see
the proof of Theorem \ref{vHopf} below). Here we used the second
base change from $\mathcal K[t]$ to $\mathcal K[t]_p^{(q)}$.
Thereby, we obtain a deformation $\mathbf u_{t,q}(\mathbf{H}
(2n;\underline1))$ of $\mathbf{H}(2n;\underline1)$ over the ring
$\mathcal{K}[t]_p^{(q)}$ with
$\dim_{\mathcal{K}}\mathbf{u}_{t,q}(\mathbf{H}(2n;\underline{1}))
=p\cdot\dim_{\mathcal{K}}\mathbf{u}(\mathbf{H}(2n;\underline{1}))
=p^{p^{2n}-1}$.

The following definition is similar to the corresponding one in \cite{HW2}.

\begin{defi}\label{fdquan}  A Hopf algebra
$(\mathbf{u}_{t,q}(\mathbf{H}(2n;\underline{1})),m,
\iota,\Delta,S,\varepsilon)$ over a ring $\mathcal K[t]_p^{(q)}$ of
characteristic $p$ is said to be a {\em finite-dimensional quantization\/}
of $\mathbf{u}(\mathbf{H}(2n;\underline{1}))$ if it is obtained from
$U(\mathbf {H}^+_\mathbb Z)[[t]]$ with its standard Hopf algebra
structure by a Drinfel'd twist, reduction modulo $p$, shrinking of
base rings, and finally a reduction modulo an appropriate ideal to ensure
``restrictedness" such that there is an isomorphism $\mathbf{u}_{t,q}
(\mathbf{H}(2n;\underline{1}))/t\mathbf{u}_{t,q}(\mathbf{H}(2n;\underline{1}))
\cong \mathbf{u}(\mathbf{H}(2n;\underline{1}))$ as algebras.
\end{defi}

Before we will be able to describe $\mathbf{u}_{t,q}(\mathbf{H}(2n;\underline{1}))$
explicitly, we need to establish another result:
\begin{lemm}\label{relation4}
Let $e:=2D_H (x^{(\ep_{-k}+2\ep_k)}),$ and $d^{(\ell)}:=\frac{1}{\ell
!}\ad e$. Then
 \begin{enumerate}
    \item[(i)] $d^{(\ell)} D_H (x^{(\alpha)})=\overline{A}_\ell D_H
    (x^{(\alpha+\ell \ep_k)})$, where $\overline{A}_\ell\in\mathbb Z_p$ as in Theorem \ref{Hopfprime}.
    \item[(ii)] $d^{(\ell)} D_H (x^{(\ep_{-i}+\ep_i)})=\delta_{\ell 0}D_H
                (x^{(\ep_{-i}+\ep_i)})-\delta_{\ell 1}\delta_{ki}e$, for $1 \leq i \leq n$.
    \item[(iii)] $d^{(\ell)}  \big(D_H (x^{(\alpha)})\big)^p=\delta_{\ell0} \big(D_H (x^{(\alpha)})\big)^p-\delta_{\ell 1}\delta_{\alpha,\ep_{-k}+\ep_k} e$.
 \end{enumerate}
\end{lemm}
\begin{proof} (i) and (ii) are easy to check.

(iii): According to Proposition 1.3(4) in Chapter 1 of \cite{HR}, the following formula holds
for any elements $a$ and $c$ in an arbitrary associative algebra with unity:
$$
c a^m=\sum\limits_{\ell=0}^{m}(-1)^\ell\binom{m}{\ell} a^{m-\ell}(\ad a)^\ell(c).
$$
Consequently, we have
\begin{equation*}
\begin{split}
 d\big(&D_Hx^{(\alpha)}\big)^p=\big[\,e,\big(D_Hx^{(\alpha)}\big)^p\,\big]\\
 &=\sum\limits_{\ell=1}^{p}(-1)^\ell
\binom{p}{\ell}\big(D_H x^{(\alpha)}\big)^{p-\ell}(\ad D_H
x^{(\alpha)})^\ell (e)\\
&\equiv (-1)^p \big(\ad D_H x^{(\alpha)}\big)^p(e)  \quad (\mod\, p)\\
&=(-1)^p \frac{1}{(\alpha !)^p}(\ad D_H x^{\alpha})^p(e)\\
 &\equiv \frac{1}{(\alpha !)^p} \prod\limits_{m=0}^{p-1} \big((m{-}2)\alpha_{-k} -
(m{-}1)\alpha_k \big)D_H
x^{p(\alpha-\ep_{-k}-\ep_k)+\ep_{-k}+2\ep_{k}} \quad (\mod\, p) \\
&\equiv \begin{cases} -e&~\text{if}~\alpha=\ep_{-k}+\ep_k\\
0&~\text{if}~\alpha \neq \ep_{-k}+\ep_k\\
\end{cases}\quad(\text{mod}\,J).
\end{split}
\end{equation*}

It follows from the definition of $d^{(\ell)}$ that $d^{(\ell)}\big((D_H
x^{(\alpha)})^p\big)=0$ for
 $2 \leq \ell \leq p-1$ in $\mathbf{u}(\mathbf{H}(2n;\underline{1}))$.
\end{proof}

As a consequence of Theorem \ref{Hopfprime}, Definition \ref{fdquan}, and
Lemma \ref{relation4}, we obtain the following result:
\begin{theorem}\label{vHopf}  For the two distinguished elements
$h:=D_H(x^{(\epsilon_{-k}+\epsilon_k)})$ and
$e:=2D_{H}(x^{(\epsilon_{-k}+2\epsilon_k)})$ $(1\leq k\leq n)$
there is a $p^{p^{2n}-1}$-dimensional noncommutative and noncocummtative Hopf algebra
$(\mathbf{u}_{t,q}(\mathbf{H}(2n;\underline{1})),m,\iota,\Delta,S,\varepsilon)$
over $\mathcal{K}[t]_p^{(q)}$ that has an undeformed algebra structure and
the following coalgebra structure resp.~antipode:
\begin{gather*}
\Delta(D_H (x^{(\alpha)})){=}D_H (x^{(\alpha)}) {\otimes}
(1{-}et)^{\alpha_k-\alpha_{-k}}{+}\sum\limits_{\ell=0}^{p-1}(-1)^\ell
h^{\lg \ell \rg} {\otimes} (1{-}et)^{-\ell}d^{(\ell)}  D_H
(x^{(\alpha)}) t^\ell, \\
S\big(
D_H(x^{(\alpha)})\big)=-(1{-}et)^{\alpha_{-k}-\alpha_{k}}\sum\limits_{\ell=0}^{p-1}
d^{(\ell)} \big(D_H(x^{(\alpha)})\big)h_{1}^{\lg \ell\rg} t^\ell,
\end{gather*}
and $\vep\big(D_H (x^{(\alpha)})\big)=0$ for $\underline{0} \le\alpha < \tau$ and $\alpha\ne\underline{0}$.
\end{theorem}
\begin{proof}
Set $U_{t,q}(\mathbf H(2n;\underline 1)):=U(\mathbf H(2n;\underline
1))\otimes_{\mathcal K}\mathcal K[t]_p^{(q)}$. Note that the statement
of Theorem \ref{Hopfprime} after the base change from $\mathcal K[t]$
to $\mathcal K[t]_p^{(q)}$ is still valid over $U_{t,q}(\mathbf H(2n;\underline 1))$.
 Denote by $I_{t,q}$ the ideal of $U_{t,q}(\mathbf H(2n;\underline
1))$ over the ring $\mathcal K[t]_p^{(q)}$ generated by the same
generators of the ideal $I$ in $U(\mathbf H(2n;\underline 1))$ after applying
a base change from $\mathcal K$ to $\mathcal K[t]_p^{(q)}$. We shall show that $I_{t,q}$ is a Hopf ideal of
$U_{t,q}(\mathbf H(2n;\underline 1))$. It suffices to verify that
$\Delta$ and $S$ preserve the generators of $I_{t,q}$ in
$U_{t,q}(\mathbf H(2n;\underline 1))$.

(I) Similar to Lemma 2.8 in \cite{HW1} or Lemma 3.4 in \cite{HW2}, we get
from Lemma~3.6:
\begin{eqnarray*}
\begin{split}
\Delta&\big((D_H x^{(\alpha)})^p \big)\\
&=\sum_{0\le j\le p \atop
\ell\ge0}\dbinom{p}{j}({-}1)^\ell(D_H x^{(\alpha)})^j h^{\langle
\ell\rangle}\otimes(1{-}et)^{j(\al_k-\al_{-k}){-}\ell}
d^{(\ell)}((D_H x^{(\alpha)})^{p{-}j})t^\ell
\end{split}
\end{eqnarray*}
\begin{eqnarray*}
\begin{split}
&= \big(D_H x^{(\alpha)}\big)^p \otimes
(1{-}et)^{p(\alpha_k-\alpha_{-k})}+
\sum\limits_{\ell=0}^{p-1}(-1)^\ell h^{\lg \ell \rg}
\otimes
(1{-}et)^{-\ell} d^{(\ell)}((D_H x^{(\alpha)})^{p})t^\ell  \\
&\equiv (D_H x^{(\al)})^p \otimes 1+1 \otimes \big(D_H x^{(\alpha)}\big)^p
+ h^{\lg 1\rg}\otimes(1{-}et)^{-1}\big(\delta_{\alpha,\ep_{-k}+\ep_k}e\big)t \quad (\mod\,p,\,I_{t,q})\\
&=(D_H x^{(\al)})^p \otimes 1+1 \otimes \big(D_H
x^{(\alpha)}\big)^p + \delta_{\alpha,\ep_{-k}+\ep_k}h \otimes (1{-}et)^{-1}et.
\end{split}
\end{eqnarray*}
So, if $\alpha \neq \ep_{-i}+\ep_i$, we get
\begin{equation*}
\begin{split}
\Delta\big((D_H x^{(\alpha)})^p \big)&=(D_H x^{(\alpha)})^p \otimes 1
+ 1 \otimes (D_H x^{(\alpha)})^p\\
& \in I_{t,q} \otimes
U_{t,q}(\mathbf{H}(2n;\underline{1}))
+U_{t,q}(\mathbf{H}(2n;\underline{1}))  \otimes I_{t,q},
\end{split}
\end{equation*}
and if $\alpha=\ep_{-i}+\ep_i$, we obtain from Theorem 3.5 and Lemma 3.8(ii):
$$
\Delta\big(D_Hx^{(\ep_{-i}+\ep_i)}\big)=D_H (x^{(\ep_{-i}+\ep_i)})
\otimes 1 +1 \otimes D_H (x^{(\ep_{-i}+\ep_i)})+\delta_{ki}h \otimes
(1{-}et)^{-1}et.
$$
Combining this with the above formula, we obtain
\begin{equation*}
\begin{split}
\Delta\big((D_H x^{(\ep_{-i}+\ep_i)})^p-D_Hx^{(\ep_{-i}+\ep_i)}\big)
&\equiv
\big((D_H x^{(\ep_{-i}+\ep_i)})^p-D_Hx^{(\ep_{-i}+\ep_i)}\big)\otimes 1 \\
&\quad +1 \otimes \big((D_H x^{(\ep_{-i}+\ep_i)})^p-D_Hx^{(\ep_{-i}+\ep_i)} \big)\\
& \in I_{t,q}\otimes U_{t,q}(\mathbf{H}(2n;\underline{1}))
+U_{t,q}(\mathbf{H}(2n;\underline{1}))\otimes I_{t,q}.
\end{split}
\end{equation*}
Thereby, we have proved that the ideal $I_{t,q}$ is a coideal of the
Hopf algebra $U_{t,q}(\mathbf{H}(2n;\underline{1}))$.

(II) Similar to Lemma 2.8 in \cite{HW1} or Lemma 3.4 in \cite{HW2}, we get
from Lemma~3.8:
\begin{eqnarray*}\label{vanti}
\begin{split}
S\big((D_H x^{(\alpha)})^p\big)&=(-1)^p
(1{-}et)^{-p(\alpha_k-\alpha_{-k})}\sum\limits_{\ell=0}^{p-1}
d^{(\ell)} (D_Hx^{(\alpha)})^p h_{1}^{\lg \ell \rg} t^\ell \\
&\equiv -\big(D_H x^{(\alpha)}\big)^p +(-1)
(-1)\delta_{\alpha,\ep_{-k}+\ep_k}e h_{1}^{\lg 1 \rg}t \\
&=-\big(D_H x^{(\alpha)}\big)^p+\delta_{\alpha,\ep_{-k}+\ep_k}e
h_{1}^{\lg 1 \rg}t.
\end{split}
\end{eqnarray*}
Therefore, when $\alpha \neq \ep_{-i}+\ep_i$, we get $S\big((D_H
x^{(\alpha)})^p\big)=-(D_H x^{(\alpha)})^p \in I_{t,q}$.

When $\alpha=\ep_{-i}+\ep_i$, by Theorem \ref{Hopfprime}, we have
\begin{equation*}
S\big(D_Hx^{(\ep_{-i}+\ep_i)}\big)=-D_H
x^{(\ep_{-i}+\ep_i)}+\delta_{ki}eh_{1}^{\lg 1 \rg}t.
\end{equation*}
Consequently, we obtain
 \begin{equation*}
S\big((D_H x^{(\ep_{-i}+\ep_i)})^p - D_Hx^{(\ep_{-i}+\ep_i)} \big)=
-\big((D_H x^{(\ep_{-i}+\ep_i)})^p - D_Hx^{(\ep_{-i}+\ep_i)} \big)
\in I_{t,q}.
\end{equation*}

Thereby, the ideal $I_{t,q}$ is indeed preserved by the antipode $S$
of the quantization $U_{t,q}(\mathbf{H}(2n;\underline{1}))$.

(III) It is obvious to see that
$\varepsilon((D_{H}x^{(\a)})^p)=0$ for all $\alpha\ne\underline{0}$ with
$\underline{0}\le\alpha<\tau$.

In other words, we have proved that $I_{t,q}$ is a Hopf ideal in
$U_{t,q}(\mathbf{H}(2n;\underline{1}))$. We thus obtain the desired
$t$-deformation of
$\mathbf{u}_{t,q}(\mathbf{H}(2n;\underline{1}))$ for the
restricted simple Hamiltonian algebra $\mathbf{H}(2n;\underline{1})$.
\end{proof}

\begin{remark} (i) \
Set $f=(1-et)^{-1}$. From Lemma \ref{relation4} and Theorem
\ref{vHopf} one gets
$$
[h,f]=f^2-f\,,\quad h^p=h\,, \quad f^p=1\,, \quad \Delta(h)=h\otimes
f+1\otimes h\,,
$$
where $f$ is a group-like element, $S(h)=-hf^{-1}$, and
$\varepsilon(h)=0$. So the subalgebra generated by $h$ and $f$ is a
Hopf subalgebra of $\mathbf{u}_{t,q}(\mathbf{H}(2n;\underline{1}))$
that is isomorphic to the well-known Radford Hopf algebra over
$\mathcal{K}$ in characteristic  $p$ (see \cite{DR}).

\smallskip
(ii) \ Given a parameter
$q\in\mathcal{K}$, one can specialize $t$ to any root of the
$p$-polynomial $t^p-qt\in\mathcal{K}[t]$ in a splitting field of
$\mathcal{K}$. For instance, by choosing $q=1$, one can specialize
$t$ to any scalar in $\mathbb{Z}_p$. For $t=0$ we get the
original standard Hopf algebra structure of
$\mathbf{u}(\mathbf{H}(2n;\underline{1}))$. In this way, we indeed
get a new Hopf algebra structure for the restricted universal
enveloping algebra $\mathbf{u}(\mathbf{H}(2n;\underline{1}))$ over
an algebraically closed field $\mathcal{K}$ consisting of the old algebra
structure and a new coalgebra structure induced
by Theorem \ref{vHopf} but of dimension $p^{p^{2n}-2}$.
\end{remark}

\subsection{More quantizations}
We can consider the modular reduction process for the quantizations of
$U(\mathbf{H}^+)[[t]]$ arising from certain products of pairwise
different and mutually commutative basic Drinfel'd twists. We will
then obtain many new families of noncommutative and noncocommutative
Hopf algebras of dimension $p^{p^{2n}-1}$ with indeterminate $t$ or
of dimension $p^{p^{2n}-2}$ by specializing $t$ to an element in
$\mathcal{K}$.

Let $A(k)_\ell$ and $A(k')_n$ denote the coefficients of the
corresponding quantizations of $U(\mathbf{H}^+_{\mathbb{Z}})$ over
$U(\mathbf{H}^+_{\mathbb{Z}})[[t]]$ given by Drinfel'd twists
$\mathcal{F}(k)$ and $\mathcal{F}(k')$, respectively (see Corollary \ref{zHopf}).
Note that $A(k)_0=A(k')_0$ $=1$, $A(k)_{-1}=A(k')_{-1}=0$.

\begin{lemm}\label{mrelation1}
For two pairs of distinguished elements $h(k):=D_H (x^{\ep_{-k}+\ep_k})$ and $e(k):=D_H
(x^{\ep_{-k}+2\ep_k})$ as well as $h(k'):=D_H(x^{\ep_{-k'}+\ep_{k'}})$ and $e(k'):=D_H
(x^{\ep_{-k'}+2\ep_{k'}})$, where $1 \leq k\ne k'\leq n$, the coalgebra structure and the antipode
of the integral quantization of $U(\mathbf{H}_{\mathbb{Z}}^{+})[[t]]$ by the
Drinfel'd twist $\mathcal{F}=\mathcal{F}(k)\mathcal{F}(k')$ with an undeformed
algebra structure are given by
\begin{eqnarray*}
\begin{split}
\Delta(D_H (x^\alpha))&=D_H (x^\alpha) \otimes (1{-}e(k)t)^{\alpha_k-\alpha_{-k}}(1{-}
e(k')t)^{\alpha_{k'}-\alpha_{-k'}} \\
&\qquad+\sum\limits_{n,\ell=0}^{\infty}  (-1)^{\ell+n}A(k')_\ell A(k)_n h(k')^{\lg \ell \rg}  h(k)^{\lg n \rg} \otimes \\
&\qquad\qquad(1{-}e(k')t)^{-\ell} (1{-}e(k)t)^{-n}
D_H(x^{\alpha+\ell \ep_{k'}+n \ep_k} t^{n+\ell}),\\
S(D_H (x^{\alpha}))&=-(1{-}e(k')t)^{\alpha_{-k'}-\alpha_{k'}}
(1{-}e(k)t)^{\alpha_{-k}-\alpha_k} \times\\
&\qquad\times\sum\limits_{n,\ell \geq 0} A(k')_\ell A(k)_n D_H (x^{\alpha+\ell \ep_k + n\ep_k}) h(k)_{1}^{\lg n \rg} h(k')_{1}^{\lg\ell\rg} t^{n+\ell},\\
\end{split}
\end{eqnarray*}
and $\varepsilon(D_H x^\alpha)=0$ for $D_H x^\alpha \in {\mathbf{H}}_{\mathbb{Z}}^{+}$.
\end{lemm}

\begin{proof}
By using Corollary \ref{zHopf}, we can get:
\begin{eqnarray*}
&&\Delta(D_H (x^\alpha))=\mathcal{F}(k)\mathcal{F}(k') \Delta_0(D_H
(x^\alpha))
\mathcal{F}(k')^{-1}\mathcal{F}(k)^{-1}\\
&&\qquad =\mathcal{F}(k) \Big( D_H (x^\alpha)\otimes
(1-e(k')t)^{\alpha_{k'}-\alpha_{-k'}}+\\
&&\qquad\quad\sum\limits_{\ell \geq 0} (-1)^\ell h(k')^{\lg \ell \rg}
\otimes (1-e(k')t)^{-\ell} A(k')_\ell D_H (x^{\alpha+\ell
 \ep_{k'}}) t^\ell \Big)\mathcal{F}(k)^{-1}\\
&&\qquad=\mathcal{F}(k)D_H (x^\alpha) \otimes
(1-e(k'))^{\alpha_{k'}-\alpha_{-k'}}F(k)\\
&&\qquad\quad +\mathcal{F}(k)\sum\limits_{\ell \geq 0} (-1)^\ell h(k')^\ell
\otimes (1-e(k')t)^{-\ell} A(k')_\ell D_H (x^{\alpha+\ell \ep_{k'}})t^\ell F(k).
\end{eqnarray*}
By Lemma \ref{vrelation2}, we can obtain
\begin{eqnarray*}
&&\mathcal{F}(k)D_H (x^\alpha) \otimes
(1-e(k'))^{\alpha_{k'}-\alpha_{-k'}}F(k)\\
&&\quad =\mathcal{F}(k)(D_H (x^\alpha) \otimes 1) F(k) \big(1 \otimes
(1-e(k')t)^{\alpha_{k'}-\alpha_{-k'}}\big)\\
&&\quad =\mathcal{F}(k)F(k)_{\alpha_{-k}-\alpha_k} \big(D_H (x^\alpha)
\otimes 1\big)(1 \otimes (1-e(k')t)^{\alpha_{k'}-\alpha_{-k'}})\\
&&\quad =D_H (x^\alpha) \otimes
(1-e(k)t)^{\alpha_k-\alpha_{-k}}(1-e(k')t)^{\alpha_{k'}-\alpha_{-k'}},\\
\mathrm{and}&\\
&&\mathcal{F}(k)\sum\limits_{\ell \geq 0} (-1)^\ell h(k')^\ell
\otimes (1-e(k')t)^{-\ell} A(k')_\ell D_H (x^{\alpha+\ell
\ep_{k'}})t^\ell F(k)\\
&&=\sum\limits_{\ell \geq 0}(-1)^\ell \big(h(k')^{\lg \ell \rg}
\otimes (1-e(k')t)^{-\ell}\big) A(k')_\ell \mathcal{F}(k) \big(1
\otimes D_H
(x^{\alpha+\ell \ep_{k'}})\big)t^{\ell} F(k)\\
&&=\sum\limits_{\ell \geq 0}(-1)^\ell h(k')^{\lg \ell \rg}
\otimes (1-e(k')t)^{-\ell}A(k')_\ell \mathcal{F}(k) \times\\
&& \qquad \sum\limits_{n
\geq 0}(-1)^n F(k)_n \big( h(k)^{\lg n \rg} \otimes d_k^{(n)} (D_H
(x^{\alpha+\ell \ep_{k'}}))\big) t^n t^{\ell} \\
&&=\sum\limits_{\ell \geq 0}(-1)^\ell \big( h(k')^{\lg \ell \rg}
\otimes (1-e(k')t)^{-\ell}\big) A(k')_\ell \mathcal{F}(k) \times\\
&&\qquad \times
\sum\limits_{n \geq 0}(-1)^n F(k)_n \big(h(k)^{\lg n \rg} \otimes
A(k)_n D_H (x^{\alpha+\ell
\ep_{k'}+n \ep_k})  t^n \big)t^{\ell} \\
&&=\sum\limits_{n,\ell \geq 0}(-1)^{\ell+n} \big(h(k')^{\lg \ell
\rg} \otimes (1-e(k')t)^{-\ell}\big) A(k')_\ell \times\\
&&\qquad  \big(1\otimes
(1-e(k)t)^{-n}\big) \big(h(k)^{\lg n \rg} \otimes A(k)_n D_H
(x^{\alpha + \ell
\ep_{k'} +n \ep_{k}}) t^{n+\ell}\big)\\
&&=\sum\limits_{n,\ell \geq 0}(-1)^{\ell+n} h(k')^{\lg \ell \rg} h(k)^{\lg n
\rg} \otimes (1-e(k')t)^{-\ell}(1-e(k)t)^{-n}\times\\
&&\quad  A(k')_\ell A(k)_n D_H
(x^{\alpha + \ell \ep_{k'}+n\ep_k}) t^{n+\ell}.
\end{eqnarray*}Hence we obtain the first statement. For the second one we have
\begin{eqnarray*}
&&S(D_H (x^\alpha))=-v(k)v(k')D_H (x^\alpha) u(k')u(k)\\
&&=-v(k) (1-e(k')t)^{\alpha_{-k'}-\alpha_{k'}}\sum\limits_{\ell \geq
0} A(k')_\ell D_H (x^{\alpha+\ell \ep_{k'}})h(k')_{1}^{\lg \ell \rg}
t^\ell u(k)\\
&&=-v(k) (1-e(k')t)^{\alpha_{-k'}-\alpha_{k'}} \sum\limits_{\ell \geq
0} A(k')_\ell u(k)_{\alpha_k-\alpha_{-k}}\times\\
&& \qquad \times \sum\limits_{n \geq 0} A(k)_n D_H (x^{\alpha + \ell \ep_{k'}
+n \ep_k})h(k)_1^{\lg n
\rg}t^n h(k')_1^{\lg \ell \rg} t^\ell \\
&&=-(1-e(k')t)^{\alpha_{-k'}-\alpha_{k'}}(1-e(k)t)^{\alpha_{-k}-\alpha_k} \cdot
\\
&&\qquad\sum\limits_{n,\ell \geq 0} A(k')_\ell A(k)_n D_H (x^{\alpha +
\ell \ep_{k'} +n \ep_k})h(k)_1^{\lg n \rg} h(k')_1^{\lg \ell \rg}
t^{n+\ell}.
\end{eqnarray*}
This completes the proof.
\end{proof}

Set $d_{k}^{(\ell)}:=\frac{1}{\ell!}(\text{\rm ad}\,e(k))^{\ell}$.
Denote the coefficients $\bar A_\ell$, $ A_\ell$ in Theorem
\ref{Hopfprime} now as $\bar A(k)_\ell$,  $ A(k)_\ell$. From Lemma
\ref{mrelation1} we get a new quantization of
$U(\mathbf{H}(2n;\underline{1}))$ over
$U_t(\mathbf{H}(2n;\underline{1}))$ by the Drinfel'd twist
$\mathcal{F}=\mathcal{F}(k)\mathcal{F}(k')$ as follows.

\begin{lemm}\label{maHopf}
For two pairs of distinguished elements $h(k):=D_H x^{(\ep_{-k}+\ep_k)}$ and $e(k):=2D_H
x^{(\ep_{-k}+2\ep_k)}$ as well as $h(k'):=D_H x^{(\ep_{-k'}+\ep_{k'})}$ and $e(k'):=2D_H
x^{(\ep_{-k'}+2\ep_{k'})}$ with $k \neq k'$, the coalgebra structure and antipode of the corresponding
quantization of $U(\mathbf{H}(2n;\underline{1}))$ over
$U_t(\mathbf{H}(2n;\underline{1}))$ with an undeformed algebra structure are given by
\begin{equation*}
\begin{split}
\Delta\big(D_H (x^{(\alpha)})\big)&=D_H (x^{(\alpha)}) \otimes
(1{-}e(k')t)^{\alpha_{k'}-\alpha_{-k'}}(1{-}e(k)t)^{\alpha_k-\alpha_{-k}}\\
&\qquad+\sum\limits_{n,\ell = 0}^{p-1}(-1)^{\ell +n}\bar{A}(k')_\ell
\bar{A}(k)_n h(k')^{\lg \ell
\rg}h(k)^{\lg n \rg} \otimes\\
&\qquad \qquad(1{-}e(k')t)^{-\ell}(1{-}e(k)t)^{-n} D_H (x^{(\alpha + \ell \ep_{k'} + n
\ep_k)})t^{n+\ell},
\end{split}
\end{equation*}
\begin{equation*}
\begin{split}
S\big(D_H
(x^{(\alpha)})\big)&=(1{-}e(k')t)^{\alpha_{-k'}-\alpha_{k'}}(1{-}e(k)t)^{\alpha_{-k}-\alpha_k}\times\\
&\qquad \times\sum\limits_{n,\ell = 0}^{p-1}\bar{A}(k')_\ell
\bar{A}(k)_n D_H (x^{(\alpha + \ell \ep_{k'} + n
\ep_k)})h(k')_1^{\lg \ell \rg} h(k)_1^{\lg n \rg}t^{n+\ell},
\end{split}
\end{equation*}
and $\varepsilon\big(D_H x^{(\alpha)}\big)=0$ for $\underline{0}\le\alpha < \tau$.
\end{lemm}
\begin{proof}From Lemma \ref{mrelation1} we obtain
\begin{gather*}
  \Delta\big(D_H (x^{(\alpha)})\big)=\frac{1}{\al!}\Delta\big(D_H (x^{\alpha})\big)\qquad\qquad \qquad\qquad\qquad\qquad\qquad\qquad\qquad\qquad\qquad\\
= D_H (x^{(\alpha)}) \otimes (1{-}e(k)t)^{\alpha_k-\alpha_{-k}}(1{-}e(k')t)^{\alpha_{k'}-\alpha_{-k'}}\qquad\qquad\qquad\qquad\qquad\qquad \\
\qquad\qquad +\frac{1}{\al!} \sum\limits_{n,\ell=0}^{\infty}  (-1)^{\ell+n}\frac{(\alpha+\ell \ep_{k'}+n \ep_k)!}{\al!} A(k')_\ell A(k)_n h(k')^{\lg \ell \rg}  h(k)^{\lg n \rg} \otimes  \qquad\qquad\qquad\qquad\qquad\qquad\qquad\\
\qquad\qquad(1{-}e(k')t)^{-\ell} (1{-}e(k)t)^{-n}
D_H(x^{(\alpha+\ell \ep_{k'}+n \ep_k)}) t^{n+\ell}\\
=D_H (x^{(\alpha)}) \otimes
(1{-}e(k')t)^{\alpha_{k'}-\alpha_{-k'}}(1{-}e(k)t)^{\alpha_k-\alpha_{-k}}\qquad\qquad\qquad\qquad\qquad\qquad\\
\qquad+\sum\limits_{n,\ell = 0}^{p-1}(-1)^{\ell +n}\bar{A}(k')_\ell
\bar{A}(k)_n h(k')^{\lg \ell
\rg}h(k)^{\lg n \rg} \otimes\qquad\qquad\qquad\qquad\qquad\qquad\\
\qquad\qquad \qquad(1{-}e(k')t)^{-\ell}(1{-}e(k)t)^{-n} D_H (x^{(\alpha + \ell \ep_{k'} + n
\ep_k)})t^{n+\ell}.\qquad\qquad\qquad\qquad\qquad\qquad\\
\end{gather*}
The other two formulas can be proved in a similar way.
\end{proof}

\begin{lemm}\label{mrelation2} For $s\ge 1$ one has
\begin{equation*}
\begin{split}
\Delta((D_{H}(x^{(\alpha)}))^s)&=\sum_{0\le j\le s\atop 0\le n, \ell\le p-1}\dbinom{s}{j}({-}1)^{n+\ell}(D_{H}(x^{(\alpha)}))^jh(k')^{\langle
\ell\rangle}h(k)^{\langle
n\rangle}\otimes \\
&\qquad
\bigl(1{-}e(k')t\bigr)^{j(\al_{k'}-\al_{-k'})-\ell}\bigl(1{-}e(k)t\bigr)^{j(\alpha_{k}-\alpha_{-k})-n}
\times
\\
&\qquad\qquad \times\,
d_{k}^{(n)}d_{k'}^{(\ell)}((D_{H}(x^{(\alpha)}))^{s{-}j})t^{n+\ell}.
\\
S((D_{H}(x^{(\alpha)}))^s)&=(-1)^s\bigl(1{-}e(k)t\bigr)^{-s(\al_k-\alpha_{-k})}
\bigl(1{-}e(k')t\bigr)^{-s(\al_{k'}-\alpha_{-k'})}\times \\
&\qquad \times
\Bigl(\sum\limits_{n,\ell=0}^{p-1}
d_{k}^{(n)}d_{k'}^{(\ell)}((D_{H}(x^{(\alpha)}))^s)h(k')_1^{\lg l \rg} h(k)_1^{\lg n \rg}t^{n+\ell}
\Bigr).
\end{split}
\end{equation*}
\end{lemm}

\begin{proof}Similar to Lemma 2.8 in \cite{HW1} or Lemma 3.4 in \cite{HW2}, we get from
Lemma~\ref{vrelation2}:
\begin{eqnarray*}
&&\Delta((D_{H}(x^{(\alpha)}))^s)=\mathcal{F}\Bigl(D_{H}(x^{(\alpha)})\otimes
1+1\otimes D_{H}(x^{(\alpha)})\Bigr)^s\mathcal{F}^{-1}\\
&&=\mathcal{F}(k)\Bigl(\sum_{0\le j\le s\atop
\ell\ge0}\dbinom{s}{j}({-}1)^\ell(D_{H}(x^{(\alpha)}))^jh(k')^{\langle
\ell\rangle}\otimes\bigl(1{-}e(k')t\bigr)^{j(\al_{k'}-\al_{-k'}){-}\ell}\times\\
&& \qquad \times
d_{k'}^{(\ell)}  (D_{H}(x^{(\alpha)}))^{s{-}j}t^\ell\Bigr)\mathcal{F}(k)^{-1}\\
&& =\mathcal{F}(k)\Bigl(\sum_{0\le j\le s\atop
\ell\ge0}\dbinom{s}{j}({-}1)^\ell\big((D_{H}(x^{(\alpha)}))^j \otimes
1\big) \cdot \big(h(k')^{\langle
\ell\rangle}\otimes\bigl(1{-}e(k')t\bigr)^{j(\al_{k'}-\al_{-k'}){-}\ell}
\big)\times\\
&&\qquad \times\big( 1\otimes
d_{k'}^{(\ell)}  (D_{H}(x^{(\alpha)}))^{s{-}j}
t^\ell\big) \Bigr)\mathcal{F}(k)^{-1}\\
&&\quad =\mathcal{F}(k) \sum_{0\le j\le s\atop
\ell\ge0}\dbinom{s}{j}({-}1)^\ell\Big((D_{H}(x^{(\alpha)}))^j \otimes
1 \Big)\times\Big( h(k')^{\langle
\ell\rangle}\otimes \bigl(1{-}e(k')t\bigr)^{j(\al_{k'}-\al_{-k'}){-}\ell}
\Big)\times\\
&&\qquad \times\Big(\sum\limits_{n=0}^{\infty} (-1)^n F(k)_n
h(k)^{\lg n \rg} \otimes d_k^{(n)}\big(
d_{k'}^{(\ell)}(D_{H}(x^{(\alpha)}))^{s{-}j}\big) t^n   t^\ell \Big)\\
&& =\mathcal{F}(k) \sum_{0\le j\le s\atop
\ell\ge0}\dbinom{s}{j}({-}1)^\ell F(k)_{n+j(\alpha_{-k}-\alpha_k)}
\Big((D_{H}(x^{(\alpha)}))^j  h(k')^{\langle
\ell\rangle}\otimes\\
&&\qquad \bigl(1{-}e(k')t\bigr)^{j(\al_{k'}-\al_{-k'}){-}\ell} \Big)
\times \Big( \sum\limits_{n=0}^{\infty} (-1)^n h(k)^{\lg n \rg}
\otimes
d_k^{(n)}\big( d_{k'}^{(\ell)} (D_{H}(x^{(\alpha)}))^{s{-}j} \big)t^{n+\ell}
\Big)\\
&&=  \sum_{0\le j\le s\atop
\ell\ge0}\dbinom{s}{j}({-}1)^{n+\ell} (D_{H}(x^{(\alpha)}))^j
h(k')^{\langle \ell\rangle} h(k)^{\langle n \rangle} \otimes\\
&&\qquad(1-e(k)t)^{j(\alpha_{k}-\alpha_{-k})-n}
(1{-}e(k')t\bigr)^{j(\al_{k'}-\al_{-k'}){-}\ell}\\
&& \qquad d_k^{(n)}\Big(
d_{k'}^{(\ell)}(D_{H}(x^{(\alpha)}))^{s{-}j} \Big)t^{n+\ell}.
\end{eqnarray*}

Similarly,
\begin{eqnarray*}
&&S((D_{H}(x^{(\alpha)}))^s)=u^{-1}S_0((D_{H} (x^{(\alpha)}))^s)\,u \\
&&=v(k)v(k')(-1)^s (D_H (x^{(\alpha)}))^s u(k')u(k)\\
&&=(-1)^s v(k)(1-e(k')t)^{s(\alpha_{-k'}-\alpha_{k'})}
\sum\limits_{\ell \geq 0} d_{k'}^{(\ell)}(D_H(x^{(\alpha)}))^s
h(k')_1^{\lg \ell \rg} t^\ell u(k)\\
&&=(-1)^s
v(k)(1-e(k')t)^{s(\alpha_{-k'}-\alpha_{k'})}\sum\limits_{\ell \geq
0}u(k)_{s(\alpha_k-\alpha_{-k})} \\
&& \qquad \times\sum\limits_{n \geq 0} d_k^{(n)}
\Big(d_{k'}^{(\ell)}(D_H(x^{(\alpha)}))^s \Big) h(k)_1^{\lg n
\rg} h(k')_1^{\lg \ell \rg} t^{n+\ell}\\
&&=(-1)^s(1-e(k)t)^{s(\alpha_{-k}-\alpha_{k})}(1-e(k')t)^{s(\alpha_{-k'}-\alpha_{k'})}\\
&& \qquad \sum\limits_{n,\ell \geq 0} d_k^{(n)}
\Big(d_{k'}^{(\ell)}(D_H(x^{(\alpha)}))^s \Big) h(k)_1^{\lg n \rg}
h(k')_1^{\lg \ell \rg} t^{n+\ell}.
\end{eqnarray*}
This completes the proof.
\end{proof}

\begin{lemm}\label{mrelation3}
Set $e(k):=2D_H(x^{(\ep_{-k}+2\ep_k)}),~e(k'):=2D_H(
x^{(\ep_{-k'}+2\ep_{k'})}),~d_k^{(n)}:=\frac{1}{n !}(\ad
e(k))^n$, and $d_{k'}^{(\ell)}:=\frac{1}{\ell !}(\ad
e(k'))^\ell$. Then

 $(\text{\rm i})$ \
 $d_k^{(n)}d_{k'}^{(\ell)}\big( D_H(x^{(\alpha)}) \big)=\bar{A}(k')_\ell
\bar{A}(k)_n D_H (x^{(\alpha+\ell \ep_{k'}+n \ep_k)})$ with $\bar{A}(k')_\ell$, $\bar{A}(k)_n\in\mathbb Z_p$.

 $(\text{\rm ii})$ \
 $d_k^{(n)}d_{k'}^{(\ell)}\big(D_H
( x^{(\ep_{-i}+\ep_i)})\big)=\delta_{\ell 0}\delta_{n0}D_H
( x^{(\ep_{-i}+\ep_i)})-\delta_{\ell
 0}\delta_{n1}\delta_{ki}e(k)-\delta_{\ell1}\delta_{k' i}\delta_{n0}e(k')$.

 $(\text{\rm iii})$
 $d_k^{(n)}d_{k'}^{(\ell)}\Big(\big(D_H (x^{(\al)})\big)^p\Big)=\delta_{\ell 0} \delta_{n0}\big(D_H (x^{(\al)})\big)^p
-\delta_{\ell 0}\delta_{n1}\delta_{\al,\ep_{-k}+\ep_k}e(k)$

\hskip4.87cm $\qquad\qquad \qquad -\,\delta_{\ell
1}\delta_{n0}\delta_{\al,\ep_{-k'}+\ep_{k'}}e(k')$.
\end{lemm}
\begin{proof} (i) For $\underline{0}\le\alpha < \tau$ we obtain from Lemma \ref{relation4}(i):
\begin{eqnarray*}
&&d_k^{(n)}d_{k'}^{(\ell)}\big(D_H (x^{(\alpha)})\big)=d_k^{(n)}
\bar{A}(k')_\ell D_H (x^{(\alpha + \ell \ep_{k'})})\\
&&=\bar{A}(k')_\ell   \frac{1}{(\alpha {+} \ell \ep_{k'})!}
d_k^{(n)}\big(D_H (x^{{\alpha + \ell \ep_{k'}}})\big)\\
&&=\bar{A}(k')_\ell \frac{1}{(\alpha {+} \ell \ep_{k'})!}
\frac{1}{n!} \prod_{j=0}^{n-1} \big((\alpha{+} \ell
\ep_{k'})_k{-}2(\alpha {+} \ell \ep_{k'})_{-k}{+}j \big) D_H
(x^{\alpha+\ell \ep_{k'}+n\ep_k})\\
&&=\bar{A}(k')_\ell \frac{1}{(\alpha {+} \ell \ep_{k'})!}
\frac{1}{n!} \prod_{j=0}^{n-1} (\alpha_k{-}2\alpha_{-k}{+}j)D_H
(x^{\alpha+\ell \ep_{k'}+n\ep_k})\\
&&=\bar{A}(k')_\ell \bar{A}(k)_n D_H (x^{(\alpha+\ell
\ep_{k'}+n\ep_k)}).
\end{eqnarray*}

(ii) By Lemma \ref{relation4}(ii), we have
\begin{equation*}
\begin{split}
d_k^{(n)}d_{k'}^{(\ell)}\big(D_H
(x^{(\ep_{-i}+\ep_i)})\big)&=d_k^{(n)}\big(\delta_{\ell 0}D_H
(x^{(\ep_{-i}+\ep_i)}) - \delta_{\ell 1}\delta_{k' i} e(k')\big)\\
\quad&=\delta_{\ell 0}\big(\delta_{n0} D_H (x^{(\ep_{-i}+\ep_i)})
-\delta_{n1}\delta_{ki}e(k)\big)-\delta_{\ell1}
\delta_{k'i}\delta_{n0}e(k')\\
\quad&=\delta_{\ell 0} \delta_{n0} D_H
(x^{(\ep_{-i}+\ep_i)})-\delta_{\ell 0} \delta_{n1}
\delta_{ki}e(k)-\delta_{\ell1} \delta_{k'i}\delta_{n0}e(k').
\end{split}
\end{equation*}

(iii) By Lemma \ref{relation4}(iii), we get
\begin{equation*}
\begin{split}
&d_k^{(n)}d_{k'}^{(\ell)}\Big(\big(D_H
(x^{(\al)})\big)^p\Big)=d_k^{(n)}\Big(\delta_{\ell 0}\big(D_H
(x^{(\alpha)})\big)^p-\delta_{\ell
1}\delta_{\alpha,\ep_{-k'}+\ep_{k'}}e(k')\Big)\\
&\qquad=\delta_{\ell 0}\Big(\delta_{n0}\big(D_H (x^{(\alpha)})\big)^p
-\delta_{n1}\delta_{\al,\ep_{-k}+\ep_k}e(k)\Big)-\delta_{\ell
1}\delta_{n0}\delta_{\alpha,\ep_{-k'}+\ep_{k'}}e(k')\\
&\qquad=\delta_{\ell 0}\delta_{n0}\big(D_H (x^{(\alpha)})\big)^p
-\delta_{\ell
0}\delta_{n1}\delta_{\al,\ep_{-k}+\ep_k}e(k)-\delta_{\ell
1}\delta_{n0}\delta_{\alpha,\ep_{-k'}+\ep_{k'}}e(k').
\end{split}
\end{equation*}

This completes the proof.
\end{proof}

By using Lemmas \ref{relation}, 3.12, \ref{mrelation2}, and \ref{mrelation3}, we obtain
a new Hopf algebra structure for the restricted universal enveloping algebra
$\mathbf{u}(\mathbf{H}(2n;\underline{1}))$ over $\mathcal{K}$ by the
product of two different and commuting vertical basic Drinfel'd twists keeping
the old algebra structure but having a new coalgebra structure and a new antipode.

\begin{theorem}
For two pairs of distinguished elements $h(k):=D_H
(x^{(\ep_{-k}+\ep_k)})$ and $e(k):=2D_H(x^{(\ep_{-k}+2\ep_k)})$ as well as
$h(k'):=D_H(x^{(\ep_{-k'}+\ep_{k'})})$ and $e(k'):=2D_H
(x^{(\ep_{-k'}+2\ep_{k'})})$ with $1 \leq k \neq k' \leq n$, there
exists a $p^{p^{2n}-1}$-dimensional noncommutative and noncocummtative Hopf algebra
$(\mathbf{u}_{t,q}(\mathbf{H}(2n;\underline{1})),m,\iota,\Delta,S,\varepsilon)$ over
$\mathcal{K}[t]_p^{(q)}$ that has an
undeformed algebra structure and the following coalgebra structure resp.~antipode:
\begin{eqnarray*}
\begin{split}
\Delta(D_H (x^{(\alpha)}))&=D_H (x^{(\alpha)}) \otimes (1{-}e(k')t)^{\alpha_{k'}-\alpha_{-k'}}(1{-}e(k)t)^{\alpha_k-\alpha_{-k}} \\
&\quad+\sum\limits_{n,\ell=0}^{p-1}  (-1)^{\ell+n} h(k')^{\lg \ell \rg}  h(k)^{\lg n \rg} \otimes
(1{-}e(k')t)^{-\ell} (1{-}e(k)t)^{-n}\times\\
&\hskip6.5cm \times d_k^{(n)} d_{k'}^{(\ell)}\big(D_H (x^{(\alpha)})\big) t^{n+\ell},\\
S(D_H (x^{(\alpha)}))&=-(1{-}e(k')t)^{\alpha_{-k'}-\alpha_{k'}}
(1{-}e(k)t)^{\alpha_{-k}-\alpha_k} \times \\
&\quad\times\sum\limits_{n,\ell = 0}^{p-1} d_k^{(n)} d_{k'}^{(\ell)}\big(D_H (x^{(\alpha)})\big) h(k')_{1}^{\lg \ell \rg} h(k)_{1}^{\lg n \rg} t^{n+\ell},
\end{split}
\end{eqnarray*}
and $\varepsilon(D_H (x^{(\alpha)}))=0$, for $\underline{0} \le\alpha<\tau$.
\end{theorem}
\begin{proof}
Let $I_{t,q}$ denote the ideal of
$(U_{t,q}(\mathbf{H}(2n;\underline{1})),m,\iota,\Delta,S,\varepsilon)$ over
the ring $\mathcal K[t]_p^{(q)}$ generated by the same generators as in $I$,
and let $q\in\mathcal{K}$.
Observe that the result in Lemma \ref{maHopf}, via the base change from
$\mathcal K[t]$ to $\mathcal K[t]_p^{(q)}$, is still valid
for $U_{t,q}(\mathbf{H}(2n;\underline{1}))$.

In what follows, we shall show that $I_{t,q}$ is a Hopf ideal of
$U_{t,q}(\mathbf{H}(2n;\underline{1}))$. To this end, it suffices to
verify that $\Delta$ and $S$ preserve the generators of $I_{t,q}$.

(I) By Lemmas \ref{mrelation2}, \ref{relation}, and \ref{mrelation3}, we obtain
\begin{eqnarray*}\label{dmrelation4}
\begin{split}
\Delta&\big(D_H (x^{(\alpha)}))^p \big)=(D_H (x^{(\alpha)}))^p \otimes
(1{-}e(k')t)^{p(\alpha_{k'}-\alpha_{-k'})}(1{-}e(k)t)^{p(\alpha_k-\alpha_{-k})} \\
&\qquad\qquad \qquad
+\sum\limits_{n,\ell=0}^{p-1}(-1)^{n+\ell}h(k')^{\lg \ell \rg}
h(k)^{\lg n \rg} \otimes (1{-}e(k')t)^{-\ell}(1{-}e(k)t)^{-n}\times \\
&\hskip5cm\times d_k^{(n)} d_{k'}^{(\ell)}\big((D_H (x^{(\alpha)}))^p\big) t^{n+\ell}\quad (\mod\,p, I_{t,q})\\
&\equiv (D_H (x^{(\alpha)}))^p  \otimes1+\sum\limits_{n,\ell=0}^{p-1} (-1)^{n+\ell} h(k')^{\lg \ell \rg}
h(k)^{\lg n \rg} \otimes (1{-}e(k')t)^{-\ell}(1{-}e(k)t)^{-n}\times\\
&\qquad\quad\times \big(\delta_{\ell 0} \delta_{n0}\big(D_H (x^{(\al)})\big)^p {-}\delta_{\ell0}\delta_{n1}\delta_{\al,\ep_{-k}{+}\ep_k}e(k){-}\delta_{\ell
1}\delta_{n0}\delta_{\al,\ep_{-k'}{+}\ep_{k'}}e(k')\big)t^{n+\ell}\\
&= (D_H (x^{(\alpha)}))^p\otimes 1 + 1 \otimes (D_H (x^{(\alpha)}))^p+\delta_{\alpha,\ep_{-k}+\ep_k} h(k) \otimes (1{-}e(k)t)^{-1}e(k)t\\
&\hskip5.35cm +\delta_{\alpha,\ep_{-k'}+\ep_{k'}} h(k')\otimes (1{-}e(k')t)^{-1}e(k')t.
\end{split}
\end{eqnarray*}

Hence, when $\alpha \neq\ep_{-i}+\ep_i$, we get
\begin{equation*}
\begin{split}
\Delta\big( (D_H (x^{(\alpha)}))^p \big)&\equiv (D_H( x^{(\alpha)}))^p
\otimes
1+1 \otimes (D_H (x^{(\alpha)}))^p\\
&\in  I_{t,q}\otimes
U_{t,q}(\mathbf{H}(2n;\underline{1}))+U_{t,q}(\mathbf{H}(2n;\underline{1}))\otimes
I_{t,q}.
\end{split}
\end{equation*}
When $\alpha=\ep_{-i}+\ep_i$,  by Lemma \ref{mrelation2} for $s=1$ and Lemma \ref{mrelation3}(ii), we have
\begin{equation*}
\begin{split}
&\Delta(D_H (x^{(\ep_{-i}+\ep_i)}))=D_H (x^{(\ep_{-i}+\ep_i)}) \otimes 1
 +\sum\limits_{n,\ell=0}^{p-1}
(-1)^{n+\ell} h(k')^{\lg \ell \rg}
h(k)^{\lg n \rg} \otimes (1{-}e(k')t)^{-\ell}\times\\
&\qquad \times(1{-}e(k)t)^{-n} \big(\delta_{\ell
0}\delta_{n0}D_H ( x^{(\ep_{-i}+\ep_i)})-\delta_{\ell
 0}\delta_{n1}\delta_{ki}e(k)-\delta_{1\ell}\delta_{k'i}\delta_{n0}e(k')\big)t^{n+\ell}\\
&=D_H (x^{(\ep_{-i}+\ep_i)}) \otimes 1 +1 \otimes D_H (x^{(\ep_{-i}+\ep_i)})+\delta_{ki}h(k)\otimes(1{-}e(k)t)^{-1}e(k)t\\
&\hskip5.84cm+\delta_{k'i}h(k')\otimes
(1{-}e(k')t)^{-1} e(k')t.
\end{split}
\end{equation*}

By combining this with the above formula, we obtain
\begin{equation*}
\begin{split}
\Delta\big((D_H (x^{(\ep_{-i}+\ep_i)}))^p-D_H (x^{(\ep_{-i}+\ep_i)})
\big)&=\big((D_H (x^{(\ep_{-i}+\ep_i)}))^p-D_H (x^{(\ep_{-i}+\ep_i)})
\big) \otimes 1\\
&\quad+ 1 \otimes \big((D_H (x^{(\ep_{-i}+\ep_i)}))^p-D_H
(x^{(\ep_{-i}+\ep_i)}) \big)\\
&\in I_{t,q}\otimes
U_{t,q}(\mathbf{H}(2n;\underline{1}))+U_{t,q}(\mathbf{H}(2n;\underline{1}))\otimes
I_{t,q}.
\end{split}
\end{equation*}
Thereby, we have proved that the ideal $I_{t,q}$ is a coideal of the
Hopf algebra $U_{t,q}(\mathbf{H}(2n;\underline{1}))$.

(II) By Lemmas \ref{mrelation2}, \ref{relation} and
\ref{mrelation3}, we have
\begin{eqnarray*}\label{amrelation}
\begin{split}
S\big((D_H (x^{(\alpha)}))^p \big)&=(-1)^p
(1{-}e(k')t)^{p(\alpha_{-k'}-\alpha_{k'})}(1{-}e(k)t)^{p(\alpha_{-k}-\alpha_{k})}\times\\
&\qquad \times\sum\limits_{n,\ell=0}^{p-1}d_k^{(n)}
d_{k'}^{(\ell)}\big((D_H (x^{(\alpha)}))^p\big)h(k')_1^{\lg \ell \rg}
h(k)_1^{\lg n \rg}  t^{n+\ell}\quad (\mod\, p, I_{t,q})\\
&\equiv-(D_H(x^{(\alpha)}))^p+\delta_{\al,\ep_{-k}+\ep_k}e(k)h(k)_{1}^{\lg 1 \rg
}t+\delta_{\al,\ep_{-k'}+\ep_{k'}}e(k')h(k')_1^{\lg 1 \rg}t.
\end{split}
\end{eqnarray*}
Hence, when $\alpha \neq\ep_{-i}+\ep_i$, we get
\begin{equation*}
\begin{split}
S\big((D_H (x^{(\alpha)}))^p\big)=-(D_H (x^{(\alpha)}))^p  \in I_{t,q}.
\end{split}
\end{equation*}

When $\alpha=\ep_{-i}+\ep_i$, by Lemma \ref{mrelation2} for $s=1$ and Lemma \ref{mrelation3}(ii), we have
\begin{equation*}
\begin{split}
S&(D_H (x^{(\ep_{-i}+\ep_i)}))=-\sum\limits_{n,\ell=0}^{p-1} d_k^{(n)}
d_{k'}^{(\ell)}\big(D_H (x^{(\ep_{-i}+\ep_i)})\big)h(k')_1^{\lg
\ell \rg} h(k)_1^{\lg n \rg}  t^{n+\ell}\\
&={-}\sum\limits_{n,\ell=0}^{p-1} \big(\delta_{\ell
0}\delta_{n0}D_H
 (x^{(\ep_{-i}+\ep_i)}){-}\delta_{\ell
 0}\delta_{n1}\delta_{ki}e(k){-}\delta_{\ell1}\delta_{n0}\delta_{k' i}e(k')\big)h(k')_1^{\lg
\ell \rg} h(k)_1^{\lg n \rg}  t^{n+\ell}\\
&=-D_H (x^{(\ep_{-i}+\ep_i)})+\delta_{ki}e(k)h(k)_1^{\lg 1
\rg}t+\delta_{k'i}e(k')h(k')_1^{\lg 1 \rg}t.
\end{split}
\end{equation*}

By combining this with the above formula, we obtain
$$
S\big((D_H (x^{(\ep_{-i}+\ep_i)}))^p-D_H
(x^{(\ep_{-i}+\ep_i)})\big)=-\big((D_H (x^{(\ep_{-i}+\ep_i)}))^p-D_H
(x^{(\ep_{-i}+\ep_i)})\big)~\in I_{t,q}.$$

Thereby, we show that the ideal $I_{t,q}$ is indeed preserved by the
antipode $H$ of the quantization
$U_{t,q}(\mathbf{H}(2n;\underline{1}))$.

\smallskip
(III) It is obvious to notice that
$\varepsilon((D_{H}(x^{(\al)}))^p)=0$ for all $\underline{0}\le\alpha<\tau$.

\smallskip
This completes the proof.
\end{proof}

\section{Quantizations of horizontal type for Lie bialgebras of Cartan type $H$}

\subsection{Quantizations of horizontal type of $\mathbf u(\mathbf{H}(2n;\underline{1}))$.}
In this section we assume that $n\geq 2$. Consider $h:=D_H
(x^{\ep_{-k}+\ep_k})$ and $e:=D_H (x^{\epsilon_k+\epsilon_m})$, where $1\leq
k, |m| \leq n, m \neq \pm k$, and denote by $\mathcal{F}(k; m)$ the
corresponding horizontal basic Drinfel'd twist. Set $d^{(\ell)}:=\frac{1}{\ell!} (\ad e)^\ell$.
For $m \in \{-n,\ldots,-1,1,\ldots, n\}$ set $\sigma(m):=1$ for $m <0$ and $\sigma(m):=-1$ for $m>0$.
By using the horizontal Drinfel'd twists, we will obtain some new
quantizations of horizontal type for the restricted universal enveloping
algebra of the Hamiltonian  algebra $\mathbf{H}(2n;\underline{1})$.
The twisted structures given by the twists $\mathcal{F}(k;m)$ on
the subalgebra $\mathbf{H}(2n;\underline{1})_0$ are the same as those on
the sympletic Lie algebra $\mathfrak{sp}_{2n}$ over a field
$\mathcal{K}$ with $\text{char}(\mathcal{K})=p$ obtained by using the
Jordanian twists $\mathcal{F}:=\mathrm{exp}(h\otimes \sigma)$, where
$\sigma:=\mathrm{ln}(1{-}e)$ for some two-dimensional carrier
subalgebra $B(2)=\text{Span}_{\mathcal K}\{h, e\}$ discussed in
by Kulish et al (see \cite{KL}, \cite{KLS}, etc.)

\begin{lemm}\label{hrelation1} For $h:=D_H (x^{\ep_{-k}+\ep_k})$ and $e:=D_H
(x^{\epsilon_k+\epsilon_m})$ $(1\leq k,|m| \leq n,\,  m \neq \pm
k)$, and $a\in \mathbb{F}, D_H
(x^\al),\ a_i \in \mathbf{H}$, the following identities hold
in $U(\mathbf{H}):$
\begin{gather*}
D_H (x^{\al}){\cdot} h_a^{\lg s \rg}= h_{a+(\al_{-k}-\al_k)}^{\lg s \rg} {\cdot} D_H( x^{\al}), \quad D_H (x^\al) {\cdot} h_{a}^{[s]}=h_{a+(\al_{-k}-\al_k)}^{[s]} {\cdot} D_H
(x^\al),\tag{\text{\rm i}}\\
 d^{(\ell)}( D_H (x^\al))=\sum\limits_{j=0}^{\ell} A_j
B_{\ell-j} D_H
(x^{\al+(\ell-j)(\ep_k-\ep_{-m})+j(\ep_m-\ep_{-k})}), \tag{\text{\rm ii}}\\
d^{(\ell)}(a_1 \cdots
a_s)=\sum\limits_{\ell_1+\cdots+\ell_s=\ell}d^{(\ell_1)}(a_1)\cdots
d^{(\ell_s)}(a_s),\tag{\text{\rm iii}} \\
D_H(x^\al) \cdot  e^s=\sum\limits_{\ell=0}^{s}(-1)^\ell
\ell!\binom{s}{\ell}e^{s-\ell}\cdot  d^{(\ell)}(D_H (x^\al)),
\tag{\text{\rm
iv}}
\end{gather*}
where $A_j=\frac{(-1)^j}{j!}\prod\limits_{i=0}^{j-1}(\al_{-k}{-}i)\in\mathbb Z$, $B_j=\frac{\sigma(m)^{j}}{j!}
\prod\limits_{i=0}^{j-1}(\al_{-m}{-}i)\in\mathbb Z$ with $A_0=B_0=1$, $A_j=0$, for $j>\al_{-k}$, $B_{j}=0$ for $j>
\al_{-m}$.
\end{lemm}
\begin{proof} We only prove (ii) as the proof of the other identities is the same as
 in \cite{HW2}.

(ii) Use induction on $\ell$. This holds for $\ell=1$, since
\begin{equation*}
\begin{split}
d(D_H( x^\al))&=[D_H (x^{\ep_k+\ep_m}), D_H (x^\al)]\\
&=\sum\limits_{i=1}^{n}\big(\pa_{-i}(\ep_k{+}\ep_{m})\al_i-\pa_{i}(\ep_k{+}\ep_{m})\al_{-i}\big)D_H
(x^{\al+\ep_k+\ep_m-\ep_{-i}-\ep_i})\\
&=\begin{cases} -\al_{-m}D_H (x^{\al+\ep_k-\ep_{-m}})-\al_{-k}D_H (x^{\al+\ep_m-\ep_{-k}}),&m >0\\
                \al_{-m}D_H (x^{\al+\ep_k-\ep_{-m}})-\al_{-k}D_H (x^{\al+\ep_m-\ep_{-k}}),&m<0
                \end{cases}\\
&=\sigma(m)\al_{-m}D_H (x^{\al+\ep_k-\ep_{-m}})-\al_{-k}D_H
(x^{\al+\ep_m-\ep_{-k}}).
\end{split}
\end{equation*}
For $\ell \geq 1$ we have
{\setlength{\arraycolsep}{0pt}
\begin{eqnarray*}
&&d^{(\ell{+}1)} D_H (x^\al)= \frac{1}{\ell{+}1}
\sum\limits_{j=0}^{\ell} A_j B_{\ell{-}j} d\Big(D_H
(x^{\al{+}(\ell-j)(\ep_k{-}\ep_{-m}){+}j(\ep_m{-}\ep_{-k})})\Big)\\
&&\quad=\frac{1}{\ell{+}1}\sum\limits_{j=0}^{\ell} A_j
B_{\ell{-}j} \Big(\sigma(m)\big(\al_{-m}{-}(\ell{-}j)\big)D_H
(x^{\al{+}(\ell{-}j{+}1)(\ep_k{-}\ep_{-m}){+}j(\ep_m{-}\ep_{-k})})\\
&&\qquad - (\al_{-k}{-}j)D_H
(x^{\al+(\ell-j)(\ep_k-\ep_{-m})+(j+1)(\ep_m-\ep_{-k})})\Big)\\
&&\quad=\sum\limits_{j=0}^{\ell}
\frac{\ell{-}j{+}1}{\ell{+}1}A_j B_{\ell-j+1}D_H
(x^{\al+(\ell-j+1)(\ep_k-\ep_{-m})+j(\ep_m-\ep_{-k})})\\
&&\qquad+ \sum\limits_{j=1}^{\ell+1}\frac{j}{\ell{+}1}
A_{j} B_{\ell-j+1} D_H
(x^{\al+(\ell-j+1)(\ep_k-\ep_{-m})+j(\ep_m-\ep_{-k})})\\
&&\quad
=B_{\ell+1}D_H (x^{\al+(\ell+1)(\ep_k-\ep_{-m})})\\
&&\qquad+\sum\limits_{j=1}^{\ell}
A_jB_{\ell-j+1}\Big(\frac{\ell{-}j{+}1}{\ell{+}1}+\frac{j}{\ell{+}1}\Big)D_H
(x^{\al+(\ell-j+1)(\ep_k-\ep_{-m})+j(\ep_m-\ep_{-k})})\\
&&\qquad+A_{\ell+1}D_H
(x^{\al+(\ell+1)(\ep_m-\ep_{-k})})\\
&&\quad=\sum\limits_{j=0}^{\ell+1} A_j B_{\ell+1-j}D_H
(x^{\al+(\ell+1-j)(\ep_k-\ep_{-m})+j(\ep_m-\ep_{-k})}).
\end{eqnarray*}
}
This completes the proof.
\end{proof}

\begin{lemm}\label{hrelation2} For $h:=D_H (x^{\ep_{-k}+\ep_k})$, $e:=D_H
(x^{\epsilon_k+\epsilon_m})$, $(1\leq k, |m| \leq n$, $m \neq \pm
k)$, and $a\in \mathbb{F}, D_H
(x^\al)\in \mathbf{H}$, the following identities hold
in $U(\mathbf{H})$$:$
\begin{gather*}
 (\ad D_H (x^\al))^s(e)=\sum\limits_{i=0}^{s}(-\sigma(m))^i\binom{s}{i}
A(s{-}i{-}1,k)A(i{-}1,m)\tag{\text{\rm
i}}\times \\
\qquad \qquad \qquad \qquad \qquad \times D_H
(x^{s\al-i(\ep_{-m}+\ep_m)-(s-i)(\ep_{-k}+\ep_k)+\ep_k+\ep_{m}}),\\
\big((D_H (x^\alpha))^s \otimes 1 \big) \cdot F_a=F_{a+s(\alpha_{-k}
-\alpha_k)} \cdot \big((D_H (x^\alpha))^s \otimes
1\big),\tag{\text{\rm
ii}}\\
(D_H (x^{\alpha}))^s \cdot u_a= u_{a+s(\alpha_k -\alpha_{-k})}
\Big(\sum\limits_{\ell=0}^{\infty}d^{(\ell)}((D_H( x^{\alpha}))^s) \cdot
h_{1-a}^{\lg \ell \rg} t^\ell \Big),\tag{\text{\rm
iii}}\\
\big (1 \otimes (D_H (x^{\alpha}))^s \big) \cdot
F_a=\sum\limits_{\ell=0}^{\infty}(-1)^\ell F_{a+\ell} \cdot
\big(h_{a}^{\lg \ell \rg} \otimes d^{(\ell)}((D_H (x^\alpha))^s) t^\ell
\big),\tag{\text{\rm iv}}
\end{gather*}
where $A(i,k)=\prod\limits_{j=0}^{i} (j\al_k{-}(j{-}1)\al_{-k})$ with $A(i,k)=1$  for $i<0$.
\end{lemm}
\begin{proof} (i) Use induction on $s$. This is true for $s=1$, since
\begin{equation*}
\begin{split}
\ad D_H (x^\al)(e)&=[D_H (x^\al),~ D_H( x^{\ep_k+\ep_m})]\\
&=\al_{-k}D_H
(x^{\al+\ep_m-\ep_{-k}})-\sigma(m)\al_{-m}D_H (x^{\al+\ep_k-\ep_{-m}}).
\end{split}
\end{equation*}

For $s\geq 1$ we have
\begin{equation*}
\begin{split}
(\ad D_H & (x^\al))^{s{+}1}(e)=\sum\limits_{i=0}^{s}(-\sigma(m))^{i}
\binom{s}{i}A(s{-}i{-}1,k)A(i{-}1,m)\times \\
&\qquad\qquad \qquad \times [D_H( x^\al) ,D_H
(x^{s\al-i(\ep_{-m}+\ep_m)-(s{-}i)(\ep_{-k}+\ep_k)+\ep_{k}+\ep_{m}})],
\end{split}
\end{equation*}
where
\begin{equation*}
\begin{split}
& [D_H (x^\al), D_H
(x^{s\al-i(\ep_{-m}+\ep_m)-(s-i)(\ep_{-k}+\ep_k)+\ep_k+\ep_m })] \\
&=\big((s{-}i)\al_k{-}(s{-}i{-}1)\al_{-k}\big)D_H (x^{(s+1)\al-i(\ep_{-m}{+}\ep_m)-(s{-}i)(\ep_{-k}{+}\ep_k){-}\ep_{-k}{+}\ep_{m}})\\
&\qquad -\sigma(m)\big(i\al_m{-}(i{-}1)\al_{-m}\big)D_H
(x^{(s+1)\al-i(\ep_{-m}+\ep_m)-(s{-}i)(\ep_{-k}{+}\ep_k){+}\ep_{k}{-}\ep_{-m}}).
\end{split}
\end{equation*}
So we can get
{\setlength{\arraycolsep}{0pt}
\begin{eqnarray*}
&&(\ad D_H (x^\al))^{s{+}1}(e)=\sum\limits_{i=0}^{s}(-\sigma(m))^{i}
\binom{s}{i}A(s{-}i{-}1,k)A(i{-}1,m)\times\\
&&\ \times\Big(\big((s{-}i)\al_k{-}(s{-}i{-}1)\al_{-k}\big)D_H (x^{(s{+}1)\al{-}i(\ep_{-m}{+}\ep_m){-}(s{-}i)(\ep_{-k}{+}\ep_k){-}\ep_{-k}{+}\ep_{m}})\\
&&\ {-}\sigma(m)\big(i\al_m{-}(i{-}1)\al_{-m}\big)D_H
(x^{(s{+}1)\al{-}i(\ep_{-m}{+}\ep_m){-}(s{-}i)(\ep_{-k}{+}\ep_k){+}\ep_{k}{-}\ep_{-m}})\Big)\\
&&=\sum\limits_{i=0}^{s}(-\sigma(m))^{i}
\binom{s}{i}A(s{-}i,k)A(i{-}1,m)D_H (x^{(s{+}1)\al{-}i(\ep_{-m}{+}\ep_m){-}(s{-}i)(\ep_{-k}{+}\ep_k){-}\ep_{-k}{+}\ep_{m}})\\
&&\ {+}\sum\limits_{i=0}^{s}(-\sigma(m))^{i+1}
\binom{s}{i}A(s{-}i{-}1,k)A(i,m)D_H (x^{(s{+}1)\al{-}i(\ep_{-m}{+}\ep_m){-}(s{-}i)(\ep_{-k}{+}\ep_k){+}\ep_{k}{-}\ep_{-m}})\\
&&=\sum\limits_{i=0}^{s}(-\sigma(m))^{i}
\binom{s}{i}A(s{-}i,k)A(i{-}1,m)D_H (x^{(s{+}1)\al{-}i(\ep_{-m}{+}\ep_m){-}(s{-}i)(\ep_{-k}{+}\ep_k){-}\ep_{-k}{+}\ep_{m}})\\
&&\ {+}\sum\limits_{i=1}^{s+1}(-\sigma(m))^{i}
\binom{s}{i{-}1}A(s{-}i,k)A(i{-}1,m)D_H (x^{(s{+}1)\al{-}i(\ep_{-m}{+}\ep_m){-}(s{-}i)(\ep_{-k}{+}\ep_k){-}\ep_{-k}{+}\ep_{m}})\\
&&=\sum\limits_{i=0}^{s+1}(-\sigma(m))^i\binom{s{+}1}{i}
A(s{-}i,k)A(i{-}1,m)D_H
(x^{(s{+}1)\al{-}i(\ep_{-m}{+}\ep_m){-}(s{+}1{-}i)(\ep_{-k}{+}\ep_k){+}\ep_k{+}\ep_{m}}).
\end{eqnarray*}}

(ii) By Lemma \ref{hrelation1}, we have
\begin{equation*}
\begin{split}
\big((D_H (x^\al))^s \otimes 1 \big) \cdot F_a&=\big((D_H( x^\al))^s
\otimes 1 \big) \cdot \sum\limits_{r=0}^{\infty} \frac{1}{r!}
h_{a}^{\lg r \rg}
\otimes e^r t^r \\
&=\sum\limits_{r=0}^{\infty} \frac{1}{r!}
h_{a+s(\al_{-k}-\al_k)}^{\lg r \rg}(D_H (x^\al))^s \otimes e^r t^r\\
&=F_{a+s(\al_{-k}-\al_k)} \cdot \big( (D_H (x^\al))^s \otimes 1 \big).
\end{split}
\end{equation*}

(iii) Use induction on $s$. For $s=1$ we have
{\setlength{\arraycolsep}{0pt}
\begin{eqnarray*}
&&D_H(x^\al) \cdot  u_a=D_H (x^\al)
\big(\sum\limits_{r=0}^{\infty}\frac{(-1)^r}{r!}h_{-a}^{[r]}e^r t^r\big)\\
&&=\sum\limits_{r=0}^{\infty}\frac{(-1)^r}{r!}
h_{-a+(\al_{-k}-\al_k)}^{[r]}D_H (x^\al) \cdot e^r t^r\\
&&=\sum\limits_{r=0}^{\infty}\frac{(-1)^r}{r!}
h_{-a+(\al_{-k}-\al_k)}^{[r]} \sum\limits_{\ell=0}^{r}(-1)^\ell\ell!
 \binom{r}{\ell}e^{r-\ell}d^{(\ell)}(D_Hx^{(\al)})t^r \\
&&=\sum\limits_{r,\ell=0}^{\infty}
\frac{(-1)^{r{+}\ell}}{(r{+}\ell)!}h_{-a+(\al_{-k}-\al_k)}^{[r+\ell]}
(-1)^\ell\ell! \binom{r{+}\ell}{\ell}e^rd^{(\ell)}(D_Hx^{(\al)}) t^{r+\ell}\\
&&=\sum\limits_{r,\ell=0}^{\infty}
\frac{(-1)^r}{r!}h_{-a+(\al_{-k}+\al_k)}^{[r]}
h_{-a+(\al_{-k}-\al_k)-r}^{[\ell]}e^rd^{(\ell)}(D_Hx^{(\al)}) t^{r+\ell}\\
&&=\sum\limits_{r,\ell=0}^{\infty} \frac{(-1)^r}{r!}
h_{-a+(\al_{-k}-\al_k)}^{[r]} e^r t^r
h_{-a+(\al_{-k}-\al_k)}^{[\ell]} d^{(\ell)}(D_Hx^{(\al)})t^\ell\\
&&=u_{a+(\al_k-\al_{-k})}\cdot\sum\limits_{\ell=0}^{\infty}
 h_{-a+(\al_{-k}-\al_k)}^{[\ell]}d^{(\ell)}(D_Hx^{(\al)})t^\ell\\
&&=u_{a+(\al_k-\al_{-k})}\cdot\sum\limits_{\ell=0}^{\infty}d^{(\ell)}(D_Hx^{(\al)})
h_{-a+\ell}^{[\ell]}t^\ell\\
&&=u_{a+(\al_k-\al_{-k})}\cdot\sum\limits_{\ell=0}^{\infty} d^{(\ell)}(D_H( x^{\al}))h_{1-a}^{\lg \ell \rg}t^\ell.\\
\end{eqnarray*}}

Suppose the identity holds for $s \geq 1$. Then we have
\begin{equation*}
\begin{split}
(&D_H(x^\al))^{s+1} \cdot u_a=D_H (x^\al) \,{\cdot}\, u_{a+s(\al_k-\al_{-k})}
\sum\limits_{r=0}^{\infty} d^{(\ell)}((D_H( x^\al))^s)
h_{1-a}^{\lg \ell \rg}t^\ell\\
&=u_{a+(s+1)(\al_k-\al_{-k})}\sum\limits_{\ell'=0}^{\infty}
d^{(\ell')}(D_H (x^\al)) h_{1-a-s(\al_k-\al_{-k})}^{\lg \ell'
\rg} t^{\ell'} {\cdot} \sum\limits_{\ell=0}^{\infty} d^{(\ell)}
((D_H (x^\al))^s) h_{1-a}^{\lg \ell \rg} t^{\ell}\\
&=u_{a+(s+1)(\al_k-\al_{-k})}\sum\limits_{\ell',\ell=0}^{\infty}
d^{(\ell')}(D_H (x^\al)) d^{(\ell)}((D_H (x^\al))^s)
h_{1-a+\ell}^{\lg \ell' \rg} h_{1-a}^{\lg \ell \rg}
t^{\ell+\ell'}\\
&=u_{a+(s+1)(\al_k-\al_{-k})}\sum\limits_{\ell=0}^{\infty}
d^{(\ell)}((D_H (x^\al))^{s+1}) h_{1-a}^{\lg \ell \rg} t^{\ell}.
\end{split}
\end{equation*}

(iv) Use induction on $s$. For $s=1$ we have
\begin{equation*}
\begin{split}
(1 \otimes D_H (x^\al)) \cdot F_a&=(1 \otimes D_H (x^\al)) \cdot \sum\limits_{r=0}^{\infty}\frac{1}{r!}h_{a}^{\lg r \rg} \otimes e^rt^r\\
&=\sum\limits_{r=0}^{\infty} \frac{1}{r!} h_{a}^{\lg r \rg}
\otimes \sum\limits_{\ell=0}^{r}(-1)^\ell\ell!\binom{r}{\ell}e^{r-\ell}
d^{(\ell)}(D_H(x^{\al})) t^r\\
&=\sum\limits_{r,\ell=0}^{\infty} \frac{\ell!}{(r+\ell)!}
h_{a}^{\lg r+\ell \rg} \otimes (-1)^\ell \binom{r+\ell}{\ell}e^r
d^{(\ell)}(D_H(x^{\al}))t^{r+\ell}\\
&=\sum\limits_{r,\ell=0}^{\infty}\frac{1}{r!}h_{a+\ell}^{\langle r \rangle}
h_{a}^{\langle \ell\rangle} \otimes (-1)^\ell e^r d^{(\ell)}(D_H(x^{\al}))
t^{r+\ell}\\
&=\sum\limits_{r=0}^{\infty} (-1)^\ell F_{a+\ell}
\big(h_{a}^{\lg \ell \rg} \otimes d^{(\ell)} (D_H (x^\al))\big)t^\ell.
\end{split}
\end{equation*}

Suppose the identity holds for $s \geq 1$. Then we have
\begin{equation*}
\begin{split}
&\big(1 \otimes (D_H (x^\al))^{s+1}\big)\cdot F_a=\big(1 \otimes D_H (x^\al)
\big) \Big(\sum\limits_{\ell=0}^{\infty} (-1)^\ell F_{a+\ell}\Big(
h_{a}^{\lg
\ell \rg} \otimes d^{(\ell)} (D_H (x^\al))^s \Big) t^\ell\Big)\\
&\quad=\sum\limits_{\ell,\ell'=0}^{\infty}
(-1)^{\ell'}F_{a+\ell+\ell'}\big(h_{a+\ell}^{\lg \ell' \rg} \otimes
d^{(\ell')}(D_H (x^\al)) \big) t^{\ell'}\Big((-1)^\ell h_{a}^{\lg
\ell \rg} \otimes d^{(\ell)} (D_H (x^\al))^s
\Big) t^\ell\\
&\quad=\sum\limits_{\ell,\ell'=0}^{\infty} (-1)^{\ell+\ell'}
F_{a+\ell+\ell'}\Big( h_{a+\ell}^{\lg \ell' \rg} h_{a}^{\lg \ell \rg}
\otimes d^{(\ell')}(D_H x^{(\al)})d^{(\ell)}(D_H (x^\al))^s\Big)
t^{\ell+\ell'}\\
&\quad=\sum\limits_{\ell=0}^{\infty}(-1)^{\ell} F_{a+\ell}\Big(
h_{a}^{\lg \ell \rg} \otimes d^{(\ell)} (D_H (x^\al))^{s+1} \Big)t^{\ell}.
\end{split}
\end{equation*}

This completes the proof.
\end{proof}

\begin{lemm}\label{hzHopf}
For the two distinguished elements  $h:=D_H (x^{\ep_{-k}+\ep_k})$ and
$e:=D_H (x^{\epsilon_k+\epsilon_m})$ $(1\leq k, |m| \leq n, m \neq
\pm k)$ the coalgebra structure and the antipode of the corresponding
horizontal integral quantization of
$U(\mathbf{H}^+_{\mathbb{Z}})$ over
$U(\mathbf{H}^+_{\mathbb{Z}})[[t]]$ by the Drinfel'd twist
$\mathcal{F}(k,m)$ with an undeformed algebra structure are given by
\begin{gather}
\Delta(D_H (x^\al))=D_H (x^\al) \otimes
(1{-}et)^{\al_k-\al_{-k}}+\sum\limits_{\ell=0}^{\infty}\sum\limits_{j=0}^{\ell}(-1)^\ell A_j
B_{\ell-j} h^{\lg \ell \rg} \otimes \\
  \qquad\qquad\qquad\qquad \qquad\qquad\qquad(1{-}et)^{-\ell}D_H (x^{\al+(\ell-j)(\ep_k-\ep_{-m})+j(\ep_m-\ep_{-k})})t^\ell, \nonumber\\
 S(D_H
(x^\al))=-(1{-}et)^{\al_{-k}-\al_k}\sum\limits_{\ell=0}^{\infty}\sum\limits_{j=0}^{\ell}
A_j B_{\ell-j}\times \\
 \qquad \qquad \qquad \qquad\qquad\qquad\qquad\times D_H( x^{\al+(\ell-j)(\ep_k-\ep_{-m})+j(\ep_m-\ep_{-k})})h_{1}^{\lg
\ell \rg}t^\ell, \nonumber\\
\varepsilon(D_H
(x^{\al}))=0.
\end{gather}
\end{lemm}

For later use we need to prove the following result.
\begin{lemm}\label{hrelation3}
For $s \geq 1$ one has
\begin{equation*}
\begin{split}
\Delta\big((D_H( x^\alpha))^s \big)&=\sum_{0\le j\le s\atop
\ell\ge0}\dbinom{s}{j}({-}1)^\ell(D_H (x^\alpha))^jh^{\langle
\ell\rangle}\otimes(1{-}et)^{j(\al_k-\al_{-k}){-}\ell}\times\\
&\hskip5.3cm\times d^{(\ell)} ((D_H (x^\alpha))^{s{-}j})t^\ell.\\
S((D_H (x^{\alpha}))^s)&=
(-1)^s(1{-}et)^{-s(\al_k-\al_{-k})}\cdot\sum\limits_{\ell=0}^{\infty}
d^{(\ell)}((D_H (x^{\alpha}))^s)\cdot h_1^{\lg
\ell\rg}t^\ell.
\end{split}
\end{equation*}
\end{lemm}

Firstly, we reduce the horizontal integral quantizations of
$U(\mathbf{H}^{+}_\mathbb{Z})$ from Lemma \ref{hzHopf}
modulo $p$ to obtain the
horizontal quantizations of $U(\mathbf{H}(2n;\underline{1}))$ over
$U_t(\mathbf{H}(2n;\underline{1}))$.

\begin{theorem}\label{hHopfprime}
For the two distinguished elements $h:=D_{H}( x^{(\epsilon_{-k}+\epsilon_k)})$ and
$e:=D_{H}(x^{(\epsilon_k+\epsilon_m)})$ $(1\leq |m| \neq k \leq n)$
the coalgebra structure and the antipode of the corresponding horizontal quantization
of $U(\mathbf{H}(2n;\underline{1}))$ over $U_t(\mathbf{H}(2n;\underline{1}))$ with an
undeformed algebra structure are given by
\begin{gather*}
\Delta(D_H (x^{(\al)}))
=D_H (x^{(\al)})\otimes (1{-}et)^{\al_{k}-\al_{-k}}+\sum\limits_{\ell=0}^{p-1}
\sum\limits_{j=0}^{\ell}(-1)^\ell \bar{A}_j \bar{B}_{\ell-j} h^{\lg \ell \rg} \otimes\\
\hskip5cm (1{-}et)^{-\ell}D_H(x^{(\al+(\ell-j)(\ep_k-\ep_{-m})+j(\ep_m-\ep_{-k}))})t^\ell,\nonumber\\
S(D_H
(x^{(\al)}))=-(1{-}et)^{\al_{-k}-\al_k}\sum\limits_{\ell=0}^{p-1}\sum\limits_{j=0}^{\ell}
\bar{A}_j\bar{B}_{\ell-j}\times\\
\hskip4.5cm \times D_H (x^{(\al+(\ell-j)(\ep_k-\ep_{-m})+j(\ep_m-\ep_{-k}))})h_1^{\lg
\ell\rg}t^\ell,
\end{gather*}
and $\varepsilon(D_H(x^{(\al)}))=0$, for $\underline{0}\le \al <\tau$. Moreover,
$\bar{A}_j \equiv (-1)^j\binom{\al_m+j}{j}$
$(\mathrm{mod}\, p)$ for $0 \leq j \leq \al_{-k}$, $\bar{B}_{\ell-j}\equiv
\sigma(m)^{\ell-j}\binom{\al_k+\ell-j}{\ell{-}j}$ $(\mathrm{mod}\, p)$,
for $0
\leq \ell-j \leq \al_{-m}$, and in all
other cases $\bar{A}_j=\bar{B}_{\ell-j}=0$.
\end{theorem}
\begin{proof}
Note that the elements $\frac{1}{\alpha!} D_H (x^\al)$ in
$\mathbf{H}^+_{\mathcal{K}}$ for $\underline{0}\le\alpha<\tau$ will be identified
with $D_H (x^{(\alpha)})$ in $\mathbf{H}(2n;\underline{1})$ and those
in $J_{\underline{1}}$ (see Section 3.2) with $0$. Hence, by
Lemma \ref{hzHopf}, we get
\begin{equation*}
\begin{split}
&\Delta(D_H (x^{(\al)}))=\frac{1}{\al !}\Delta(D_H (x^\al))=D_H
(x^{(\al)})\otimes(1{-}et)^{\al_k-\al_{-k}} \\
&\qquad+\sum\limits_{\ell=0}^{p-1}\sum\limits_{j=0}^{\ell}\frac{1}{\al !}(-1)^\ell A_j
B_{\ell-j} h^{\lg
\ell\rg} \otimes (1{-}et)^{-\ell}D_H (x^{\al+(\ell-j)(\ep_k-\ep_{-m})+j(\ep_m-\ep_{-k})})t^\ell\\
&\quad= D_H
(x^{(\al)})\otimes(1{-}et)^{\al_k-\al_{-k}}+\sum\limits_{\ell=0}^{p-1}
\sum\limits_{j=0}^{\ell}
\frac{\big(\al{+}(\ell{-}j)(\ep_k{-}\ep_{-m}){+}j(\ep_m{-}\ep_{-k})\big)!}{\al!}\cdot\\
&\qquad \cdot({-}1)^\ell A_j B_{\ell-j}  h^{\lg \ell\rg} \otimes
(1{-}et)^{-\ell} D_H
\big(x^{(\al+(\ell-j)(\ep_k-\ep_{-m})+j(\ep_m-\ep_{-k}))}\big)t^\ell,
\end{split}
\end{equation*}
where the coefficients below, by the definition of $A_j$ and $B_j$, become
\begin{equation*}
\begin{split}
&\frac{\big(\al{+}(\ell{-}j)(\ep_k{-}\ep_{-m}){+}j(\ep_m{-}\ep_{-k})\big)!}{\al!}A_j B_{\ell-j}\\
&\quad=A_j \frac{(\al_{-k}{-}j)!(\al_m{+}j)!}{\al_{-k}!\al_m!}\cdot B_{\ell-j} \frac{(\al_{-m}{-}(\ell{-}j))!(\al_k{+}(\ell{-}j))!}{\al_{-m}!\al_k!}\\
&\quad=(-1)^j\binom{\al_m{+}j}{j}\cdot\sigma(m)^{\ell-j}\binom{\al_k{+}\ell{-}j}{\ell{-}j}=\bar{A}_j \bar{B}_{\ell-j}.
\end{split}
\end{equation*}
Therefore, we have verified the formulas both for the coproduct and the antipode.

This completes the proof.
\end{proof}

Before we will be able to describe $\mathbf{u}_{t,q}(\mathbf{H}(2n;\underline{1}))$
explicitly, we need to establish another result:

\begin{lemm}\label{hrelation5}
Set $e:=D_H (x^{(\ep_k+\ep_m)})$ and $d^{(\ell)}:=\frac{1}{\ell !}\ad
e$. Then

\smallskip
$(\text{\rm i})$ \ $d^{(\ell)}(D_H (x^{(\al)}))=
\sum\limits_{j=0}^{\ell} \bar{A}_j \bar{B}_{\ell-j}D_H
(x^{(\al+(\ell-j)(\ep_k-\ep_{-m})+j(\ep_m-\ep_{-k}))})$, where
$\bar{A}_j$,  $\bar{B}_j\in\mathbb Z_p$, as in Theorem \ref{hHopfprime}.

\smallskip
$(\text{\rm ii})$ \ $d^{(\ell)}(D_H
(x^{(\ep_{-i}+\ep_i)}))=\delta_{\ell, 0} D_H
(x^{(\ep_{-i}+\ep_i)})+\delta_{\ell,
1}(\delta_{i,-m}{-}\delta_{i,m}{-}\delta_{i,k})e$.

\smallskip
$(\text{\rm iii})$ \ $d^{(\ell)}\big((D_H
(x^{(\al)}))^p\big)=\delta_{\ell, 0}(D_H(x^{(\al)}))^p-\delta_{\ell,
1}\big(\delta_{\al,\ep_{-k}+\ep_k}{-}\sigma(m)\delta_{\al,\ep_{-m}+\ep_m}\big)e$.
\end{lemm}
\begin{proof} (i) By Lemma \ref{hrelation1}(ii) and the proof of Theorem \ref{hHopfprime}, we have
\begin{equation*}
\begin{split}
d^{(\ell)}(D_H (x^{(\al)}))&=\frac{1}{\al !} d^{(\ell)}(D_H (x^\al))\\
&=\sum\limits_{j=0}^{\ell}
\frac{\big(\al{+}(\ell{-}j)(\ep_k{-}\ep_{-m}){+}j(\ep_m{-}\ep_{-k})\big)!}{\al
!} A_j B_{\ell-j}\times\\
&\qquad\qquad\times D_H
(x^{(\al+(\ell-j)(\ep_k{-}\ep_{-m}){+}j(\ep_m{-}\ep_{-k}))}) \\
&= \sum\limits_{j=0}^{\ell}
\bar{A}_j\bar{B}_{\ell-j}D_H
(x^{(\al+(\ell-j)(\ep_k{-}\ep_{-m}){+}j(\ep_m{-}\ep_{-k}))}).
\end{split}
\end{equation*}

(ii) $d(D_H(x^{(\ep_{-i}+\ep_i)}))=[D_H(x^{(\ep_k+\ep_m)}), D_H
(x^{(\ep_{-i}+\ep_i)})]=(\delta_{i,-m}{-}\delta_{i,m}{-}\delta_{i,k}) e$.
So if $\ell \geq 2$, then $d^{(\ell)}(D_H(x^{(\ep_{-i}+\ep_i)}))=0$. Thus, we obtain (ii).

(iii)  For $\ell=1$ we conclude from Lemma 4.2:
\begin{equation*}
\begin{split}
d((D_H x^{(\al)})^p)&=[e,(D_H
(x^{(\al)}))^p]\\
&=\sum\limits_{\ell=1}^{p}(-1)^\ell \binom{p}{\ell}(D_H
(x^{(\al)}))^{p-\ell} \cdot \big(\ad D_H
(x^{(\al)})\big)^\ell ( e)\\
& \equiv(-1)^p  \big(\ad D_H (x^{(\al)})\big)^p( e )\qquad(\mod\, p)\\
&=-\frac{1}{(\al!)^p} (\ad D_H (x^\al))^p \Big(D_H (x^{\ep_k+\ep_m})\Big)\\
&=-\frac{1}{(\al!)^p}\sum\limits_{i=0}^{p}\binom{p}{i}
(-\sigma(m))^{i}A(p{-}i{-}1,k)A(i{-}1,m)\times \\
&\qquad \times D_H (x^{p\al
-i(\ep_{-m}+\ep_m)-(p-i)(\ep_{-k}+\ep_k)+\ep_{k}+\ep_{m}})\\
&=-\frac{1}{(\al!)^p}A(p{-}1,k)D_H
(x^{p(\al-\ep_{-k}-\ep_k)+\ep_m+\ep_k})\\
&\qquad +\frac{1}{(\al!)^p}\sigma(m)^pA(p{-}1,m)D_H (x^{p(\al-\ep_{-m}-\ep_m)+\ep_m+\ep_k})\\
&\equiv-\frac{1}{(\al!)^p}\delta_{\al,\ep_{-k}+\ep_k}D_H (x^{\ep_k+\ep_m})
+\frac{1}{(\al!)^p}\delta_{\al,\ep_{-m}+\ep_m}\sigma(m)^p D_H
(x^{\ep_k+\ep_m})\quad (\mod\, J)\\
&=-\big(\delta_{\al,\ep_{-k}+\ep_k}-\sigma(m)\delta_{\al,\ep_{-m}+\ep_m}
\big)e.
\end{split}
\end{equation*}

So if $\ell \geq 2$, then $d^{(\ell)}  (D_H (x^{(\al)}))^p=0$. Thus, we obtain (iii).
\end{proof}

By using Theorem \ref{hHopfprime} and Lemma 4.6, we can now prove the following result:
\begin{theorem}\label{hHopf}
For the two distinguished elements $h:=D_H (x^{(\ep_k+\ep_{-k})})$ and $e:=D_H
(x^{(\ep_k+\ep_m)})$ $(1 \leq k \neq |m| \leq n)$ there exists a $p^{p^{2n}-1}$-dimensional
noncommutative and noncocommutative Hopf algebra $(\mathbf{u}_{t,q}(\mathbf{H}
(2n;\underline{1})),m,\iota,\Delta,S,\varepsilon)$  (of horizontal type) over
$\mathcal{K}[t]_p^{(q)}$ that has an undeformed algebra structure and the following
coalgebra structure resp.~antipode:
\begin{gather*}
\Delta(D_H (x^{(\al)})){=}D_H (x^{(\al)}) {\otimes}
(1{-}et)^{\al_k-\al_{-k}}+\sum\limits_{\ell=0}^{p-1}(-1)^\ell h^{\lg
\ell \rg}{\otimes}
(1{-}et)^{-\ell} d^{(\ell)}  (D_H (x^{(\al)}))t^\ell, \\
S(D_H (x^{(\al)}))=-(1{-}et)^{\al_{-k}-\al_k}\sum\limits_{\ell=0}^{p-1}
d^{(\ell)}(D_H (x^{(\al)}))h_{1}^{\lg \ell \rg} t^\ell,
\end{gather*}
and $\varepsilon(D_H (x^{(\al)}))=0$, for $\underline{0}\le\alpha<\tau$.
\end{theorem}
\begin{proof}
By utilizing the same argument as in the proof of Theorem
\ref{vHopf}, we shall show that the ideal $I_{t,q}$ is a Hopf ideal
of the twisted Hopf algebra
$U_{t,q}(\mathbf{H}(2n;\underline{1}))$. To this end, it suffices to verify that $\Delta$
and $S$ preserve the generators in $I_{t,q}$.

(I) By Lemma \ref{hrelation3}, Theorem \ref{hHopfprime}, and Lemma 4.6, we obtain
\begin{equation*}
\begin{split}
&\Delta\big((D_H (x^{(\al)}))^p \big)\\
&=\sum_{0\le j\le p \atop
\ell\ge0}\dbinom{p}{j}({-}1)^\ell(D_H (x^{(\alpha)}))^j h^{\langle
\ell\rangle}\otimes(1{-}et)^{j(\al_k-\al_{-k}){-}\ell}
d^{(\ell)}((D_H (x^{(\alpha)}))^{p{-}j})t^\ell\\
&= (D_H
(x^{(\al)}))^p {\otimes}
(1{-}et)^{p(\al_{k}-\al_{-k})} \\
&\qquad+ \sum\limits_{\ell=0}^{p-1}(-1)^\ell h^{\lg \ell \rg}
\otimes (1{-}et)^{-\ell} d^{(\ell)}(D_H (x^{(\al)}))^p t^\ell \qquad(\mod\, p, I_{t,q})\\
&\equiv(D_H (x^{(\al)}))^p {\otimes} 1 +1 {\otimes}(D_H (x^{(\al)}))^p +h
{\otimes} (1{-}et)^{-1}\big(\delta_{\al,\ep_{-k}+\ep_k}{-}\sigma(m)
\delta_{\al,\ep_{-m}+\ep_m}\big) et.
\end{split}
\end{equation*}

That is,
\begin{equation}
\label{hDelta}\begin{split}
\Delta\big((D_H
(x^{(\al)}))^p\big)&=(D_H (x^{(\al)}))^p
\otimes 1 + 1 \otimes(D_H (x^{(\al)}))^p\\
&\quad +h \otimes (1{-}et)^{-1}\delta_{\al,\ep_i+\ep_{-i}}\big(
\delta_{k,i}+ \delta_{m,i}-\delta_{-m,i}\big) et.
\end{split}
\end{equation}
So, when $\alpha \neq\ep_{-i}+\ep_i$, $1 \leq i \leq n$, we have
\begin{equation*}
\begin{split}
\Delta\big((D_H (x^{(\alpha)}))^p \big)&=(D_H (x^{(\alpha)}))^p \otimes 1
+ 1 \otimes (D_H (x^{(\alpha)}))^p \\
&\in I_{t,q} \otimes
U_{t,q}(\mathbf{H}(2n;\underline{1}))
+U_{t,q}(\mathbf{H}(2n;\underline{1}))  \otimes I_{t,q};
\end{split}
\end{equation*}
and when $\alpha=\ep_{-i}+\ep_i$, it follows from Theorem \ref{hHopfprime} that
\begin{equation*}
\begin{split}
\Delta(D_H (x^{(\ep_{-i}+\ep_i)}))&=D_H (x^{(\ep_{-i}+\ep_i)}) \otimes 1+1 \otimes D_H
(x^{(\ep_{-i}+\ep_i)})\\
&\quad+h \otimes (1{-}et)^{-1}(\delta_{k,i}+\delta_{m,i}-\delta_{-m,i})et.
\end{split}
\end{equation*}
Combining this with (\ref{hDelta}), we obtain
\begin{equation*}
\begin{split}
\Delta\big((D_H x^{(\ep_i+\ep_{-i})})^p-D_Hx^{(\ep_i+\ep_{-i})}\big)
&\equiv
\big((D_H (x^{(\ep_i+\ep_{-i})}))^p-D_H(x^{(\ep_i+\ep_{-i})})\big)\otimes 1 \\
&\quad +1 \otimes \big((D_H( x^{(\ep_i+\ep_{-i})}))^p-D_H(x^{(\ep_i+\ep_{-i})})\big) \\
& \in I_{t,q} \otimes U_{t,q}(\mathbf{H}(2n;\underline{1}))
+U_{t,q}(\mathbf{H}(2n;\underline{1}))  \otimes I_{t,q}.
\end{split}
\end{equation*}

So, $I_{t,q}$ is indeed a coideal of $U_{t,q}(\mathbf{H}(2n;\underline{1}))$.

(II) By Lemma \ref{hrelation3}, Lemma \ref{hrelation5}, and Theorem \ref{hHopfprime}, we have
\begin{equation*}
\begin{split}\label{hanti}
S\big((D_H( x^{(\al)}))^p\big)&=(-1)^p
(1{-}et)^{p(\al_{-k}-\al_k)}\sum\limits_{\ell=0}^{p-1}
d^{(\ell)}\big((D_H( x^{(\al)}))^p\big) h_{1}^{\lg \ell \rg} t^\ell \quad(\mod\,p, I_{t,q})\\
&\equiv-(D_H( x^{(\al)}))^p+\big(\delta_{\al,\ep_k+\ep_{-k}}
-\sigma(m)\delta_{\al,\ep_m+\ep_{-m}}\big)e
h_{1}^{\lg 1 \rg}t.
\end{split}
\end{equation*}
So, when $\alpha \neq\ep_{-i}+\ep_i$, we
have $S\big((D_H(x^{(\al)}))^p \big)=-(D_H(x^{(\al)}))^p \in I_{t,q}$.
When $\alpha =\ep_{-i}+\ep_i$, it follows from Theorem \ref{hHopfprime} or Lemmas 4.4 and 4.6 (ii) that
\begin{equation*}
\begin{split}
S(D_H (x^{(\ep_i+\ep_{-i})}))&=-\sum\limits_{\ell=0}^{p-1}d^{(\ell)}  (D_H (x^{(\al)}))h_{1}^{\lg \ell \rg}t^\ell \\&
=-D_H
(x^{(\ep_i+\ep_{-i})})+(\delta_{m,i}+\delta_{k,i}-\delta_{-m,i})e
h_{1}^{\lg 1 \rg} t.
\end{split}
\end{equation*}

So, we obtain

$S\big((D_Hx^{(\ep_i+\ep_{-i})})^p-D_H
x^{(\ep_i+\ep_{-i})}\big)=-\big((D_Hx^{(\ep_i+\ep_{-i})})^p-D_H
x^{(\ep_i+\ep_{-i})}\big) \in I_{t,q}$.

Thereby, we have shown that $I_{t,q}$ is preserved by the antipode $S$ of
$U_{t,q}(\mathbf{H}(2n;\underline{1}))$ as in Theorem \ref{hHopfprime}.

\smallskip
So, $I_{t,q}$ is a Hopf ideal in
 $U_{t,q}(\mathbf{H}(2n;\underline{1}))$, and we get a
finite-dimensional horizontal quantization for
 $\mathbf{u}_{t,q}(\mathbf{H}(2n;\underline{1}))$.
\end{proof}

\subsection{Jordanian modular quantizations of
$\mathbf{u}(\mathfrak{sp}_{2n})$}
Let $\mathbf{u}(\mathfrak{sp}_{2n})$ denote the restricted universal
enveloping algebra of $\mathfrak{sp}_{2n}$. Since Drinfel'd twists
$\mathcal{F}(k; m)$ of horizontal type act stably on the subalgebra
$U((\mathbf H_{\mathbb Z}^+)_0)[[t]]$, and so consequently on $\mathbf
u_{t,q}(\mathbf H(2n;\underline 1)_0)$, these give rise to the
Jordanian quantizations for $\mathbf u_{t,q}(\mathfrak{sp}_{2n})$.

By Lemma \ref{hrelation5}(i), we have for $r\ne s$:
\begin{equation*}
\begin{split}
d^{(\ell)}(D_H( x^{(\ep_r+\ep_s)}))&=\delta_{\ell 0} D_H
(x^{(\ep_r+\ep_s)})\\
& {+}\delta_{\ell1}\Big({\sigma}(m)(\delta_{-m,r}(\delta_{k,s}{+}1)
{+}\delta_{-m,s}(\delta_{k,r}{+}1))D_H
(x^{(\ep_k{+}\ep_r{+}\ep_s{-}\ep_{-m})})\\
& {-}\big(\delta_{-k,r}(\delta_{m,s}{+}1)
+\delta_{-k,s}(\delta_{m,r}{+}1)\big)D_H
(x^{(\ep_m{+}\ep_r{+}\ep_s{-}\ep_{-k})})\Big)\\
& {-}\delta_{\ell2} \sigma(m)
(\delta_{-m,r}\delta_{-k,s}{+}\delta_{-m,s}\delta_{-k,r})e,\\
d^{(\ell)}(D_H (x^{(2\ep_{-k})}))&=\delta_{\ell0}D_H
(x^{{2\ep_{-k}}})-\delta_{\ell1}D_H(x^{(\ep_{-k}+\ep_m)})
+\delta_{\ell2}D_H (x^{(2\ep_m)}),\\
d^{(\ell)}(D_H (x^{(2\ep_{-m})}))&=\delta_{\ell0}D_H
(x^{(2\ep_{-m})})+\delta_{\ell1}\sigma(m)D_H
(x^{(\ep_{-m}+\ep_k)})+\delta_{\ell2}D_H (x^{(2\ep_k)}),\\
d^{(\ell)}(D_H (x^{(2\ep_{r})}))&=\delta_{\ell,
0}D_H (x^{(2\ep_{r})}) \qquad \textit{ for  } \  r\neq -k, -m.
\end{split}
\end{equation*}

The next result now follows from Theorem \ref{hHopf}:
\begin{theorem}\label{0degree} For the two distinguished elements
$h:=D_H(x^{(\epsilon_{-k}+\epsilon_k)})$ and
$e:=D_H(x^{(\epsilon_k+\ep_m)})$ $(1\leq |m|\neq k \leq n)$ the coalgebra
structure and the antipode of the corresponding Jordanian quantization of
$\mathbf{u}(\mathbf{H}(2n;\underline{1})_0)\cong \mathbf
u(\mathfrak{sp}_{2n})$ over
$\mathbf{u}_{t,q}(\mathbf{H}(2n;\underline{1})_0)\cong \mathbf
u_{t,q}(\mathfrak{sp}_{2n})$ with an undeformed algebra structure are given by
\begin{equation*}
\begin{split}
\Delta(D_H&(x^{(\ep_r+\ep_s)}))=D_H (x^{(\ep_r+\ep_s)}) \otimes (1{-}et)^{\delta_{r,k}+\delta_{s,k}-\delta_{r,-k}-\delta_{s,-k}}+1
\otimes D_H (x^{(\ep_r+\ep_s)})\\
&-h {\otimes}
(1{-}et)^{-1}\Big(
\sigma(m)(\delta_{-m,r}(\delta_{k,s}{+}1)
{+}\delta_{-m,s}(\delta_{k,r}{+}1))D_H
(x^{(\ep_k+\ep_r+\ep_s-\ep_{-m})})\\
&-\big(\delta_{-k,r}(\delta_{m,s}{+}1)
{+}\delta_{-k,s}(\delta_{m,r}{+}1)\big)D_H
(x^{(\ep_m+\ep_r+\ep_s-\ep_{-k})})\Big)t\\
&-\sigma(m)\big(\delta_{r,-m}\delta_{s,-k}{+}\delta_{s,-m}\delta_{r,-k}\big)h^{\lg 2 \rg} \otimes
(1{-}et)^{-2}e
t^2, \\
\Delta(D_H&(x^{(2\ep_r)}))=D_H (x^{(2\ep_r)}) \otimes (1{-}et)^{2\delta_{k,r}-2\delta_{-k,r}}+1 \otimes D_H
(x^{(2\ep_r)})\\
&+\delta_{r,-k}\Big(h{\otimes} (1{-}et)^{-1}D_H(x^{(\ep_{-k}+\ep_m)})t\\
&+h^{\lg 2 \rg}{\otimes} (1{-}et)^{-2} D_H( x^{(2\ep_m)}) t^2\Big)\\
&-\delta_{r,-m}\Big(\sigma(m)h {\otimes} (1{-}et)^{-1} D_H
(x^{(\ep_{-m}+\ep_k)})t\\
&-h^{\lg 2 \rg} {\otimes} (1{-}et)^{-2}D_H
(x^{(2\ep_k)})t^2\Big),
\end{split}
\end{equation*}
\begin{equation*}
\begin{split}
S(D_H&
(x^{(\ep_r+\ep_s)}))=-(1{-}et)^{\delta_{r,-k}+\delta_{s,-k}-\delta_{r,k}-\delta_{s,k}}\Big[\,D_H
(x^{(\ep_r+\ep_s)})\\
&+\left(\sigma(m)(\delta_{-m,r}(\delta_{k,s}{+}1)
{+}\delta_{-m,s}(\delta_{k,r}{+}1))D_H
(x^{(\ep_k+\ep_r+\ep_s-\ep_{-m})})\right.\\
&-\left.\big(\delta_{-k,r}(\delta_{m,s}{+}1)
{+}\delta_{-k,s}(\delta_{m,r}{+}1)\big)D_H
(x^{(\ep_m+\ep_r+\ep_s-\ep_{-k})})\right) h_{1}^{\lg 1 \rg} t\\
&-\sigma(m) (\delta_{r,-m}\delta_{s,-k}{+}\delta_{s,-m}\delta_{r,-k})e
h_{1}^{\lg 2 \rg} t^2\, \Big],\\
S(D_H&(x^{(2\ep_r)}))=-(1{-}et)^{2\delta_{r,-k}-2\delta_{r,k}}\Big[\,D_H
(x^{(2\ep_r)})\\
&- \delta_{r,-k}\Big(D_H (x^{(\ep_{-k}+\ep_{m})}) h_{1}^{\lg 1 \rg} t-D_H
(x^{(2\ep_m)})h_{1}^{\lg 2 \rg} t^2\Big)\\
&+\delta_{r,-m}\Big(\sigma(m)D_H (x^{(\ep_{-m}+\ep_k)}) h_{1}^{\lg 1 \rg} t+D_H
(x^{(2\ep_k)})h_{1}^{\lg 2 \rg} t^2\Big)\,\Big],
\end{split}
\end{equation*}
and $\varepsilon(D_H (x^{(\ep_r+\ep_s)}))=\ep(D_H (x^{(2\ep_r)}))=0$, where $r\ne s$.
\end{theorem}

\begin{remark}\label{rek}
As $\mathbf{H}(2n;\underline{1})_0\cong \mathfrak{sp}_{2n}$
via the identification of $D_{H}(x^{(\epsilon_r+\epsilon_s)})$
with $\frac{\sigma(s)E_{r,-s}+\sigma(r)E_{s,-r}}{(\epsilon_r{+}\epsilon_s)!}$ for $-n\leq r, s\leq n$, 
we get a Jordanian
quantization for $\mathfrak{sp}_{2n}$, which has been discussed by
Kulish et al (cf. \cite{KL}, \cite{KLS} etc.).
\end{remark}

\begin{coro}\label{Jordanian} For the two distinguished elements
$h:=E_{k,k}-E_{-k,-k}$ and $e:=\sigma(m)E_{k,-m}-E_{m,-k}$ $(1\leq k
\neq |m|\leq n)$ the coalgebra structure and the antipode of the
corresponding Jordanian quantization of
$\mathbf{u}(\mathfrak{sp}_{2n})$ over
$\mathbf{u}_{t,q}(\mathfrak{sp}_{2n})$ with an undeformed algebra structure are given by
\begin{equation*}
\begin{split}
&\Delta(\sigma(s)E_{r,-s}{+}\sigma(r)E_{s,-r})\\
&=\big(\sigma(s)E_{r,-s}{+}\sigma(r)E_{s,-r}\big)
{\otimes}(1{-}et)^{\delta_{r,k}+\delta_{s,k}{-}\delta_{r,-k}{-}\delta_{s,-k}}{+}1 {\otimes}
\big(\sigma(s)E_{r,-s}{+}\sigma(r)E_{s,-r}\big) \\
&-h{\otimes} (1{-}et)^{-1} \Big(\sigma(m)\delta_{r,-m}\big({-}E_{s,-k}{+}\sigma(s)E_{k,-s}\big){+}\sigma(m)\delta_{s,-m}\big({-}E_{r,-k}{+}\sigma(r)E_{k,-r}\big)\\
&-\delta_{r,-k}\big(\sigma(m)E_{s,-m}{+}\sigma(s)E_{m,-s}\big){-}\delta_{s,-k}\big(\sigma(m)E_{r,-m}{+}\sigma(r)E_{m,-r}\big)\Big)t\\
&-\sigma(m)(\delta_{r,-m}\delta_{s,-k}{+}\delta_{s,-m}\delta_{r,-k}) h^{\lg 2 \rg} {\otimes}
(1{-}et)^{-2}\big(\sigma(m)E_{k,-m}{-}E_{m,-k}\big)t^2,\\
&\Delta(\sigma(r)E_{r,-r})=\sigma(r)E_{r,-r} \otimes (1{-}et)^{2\delta_{k,r}-2\delta_{-k,r}}+1 \otimes
\sigma(r)E_{r,-r}\\
& +\delta_{r,-k}\bigg(\Big(h{\otimes}
(1{-}et)^{-1}(\sigma(m)E_{-k,-m}{+}E_{m,k})\Big)t+h^{\lg 2 \rg} {\otimes}
(1{-}et)^{-2} \sigma(m)E_{m,-m}t^2 \bigg)\\
& +\delta_{r,-m}\bigg(\Big(\sigma(m)h {\otimes}
(1{-}et)^{-1}\big(\sigma(m)E_{k,m}{+}E_{-m,-k}\big)\Big)t-h^{\lg 2 \rg}
{\otimes} (1{-}et)^{-2} E_{k,-k}t^2\bigg),\\
&S(\sigma(s)E_{r,-s}{+}\sigma(r)E_{s,-r})=-(1{-}et)^{\delta_{r,-k}
+\delta_{s,-k}-\delta_{r,k}-\delta_{s,k}}\Big((\sigma(s)E_{r,-s}{+}\sigma(r)E_{s,-r})\\
&+\Big(\sigma(m)\big(\delta_{r,-m}({-}E_{s,-k}{+}\sigma(s)E_{k,-s})
 +\delta_{s,-m}({-}E_{r,-k}{+}\sigma(r)E_{k,-r})\big)\\
&-\delta_{r,-k}(\sigma(m)E_{s,-m}{+}\sigma(s)E_{m,-s})
{-}\delta_{s,-k}(\sigma(m)E_{r,-m}{+}\sigma(r)E_{m,-r})\big)h_{1}^{\lg 1 \rg}t\\
&-\sigma(m)(\delta_{r,-m}\delta_{s,-k}
{+}\delta_{s,-m}\delta_{r,-k})(\sigma(m)E_{k,-m}{-}E_{m,-k})h_{1}^{\lg
2 \rg}t^2\Big),
\end{split}
\end{equation*}
\begin{equation*}
\begin{split}
&S(\sigma(r)E_{r,-r})=-(1{-}et)^{2\delta_{r,-k}{-}2\delta_{r,k}}\sigma(r)E_{r,-r} \\
& +\delta_{r,-k}(1{-}et)^2
\Big(\big(\sigma(m)E_{-k,-m}{+}E_{m,k}\big)h_{1}^{\lg 1 \rg} t{-}\sigma(m)E_{m,-m}h_{1}^{\lg 2 \rg}t^2\Big) \\
&-\delta_{r,-m}\sigma(m)\Big(\big(\sigma(m)E_{k,m}{+}E_{-m,-k}\big)h_{1}^{\lg 1 \rg} t{+}E_{k,-k}h_{1}^{\lg 2 \rg}t^2\Big),\\
&+\delta_{r,-m}\big((E_{k,m}+\sigma(m)E_{-m,-k})h_1^{\langle 1 \rangle}t+E_{k,-k}h_1^{\langle 2 \rangle}t^2\big),\end{split}
\end{equation*}
and $\varepsilon(\sigma(s)E_{r,-s}{+}\sigma(r)E_{s,-r})=\ep(\sigma(r)E_{r,-r})=0$ for $1 \leq  |r|,|s| \leq n$ and $r \neq s$.
\end{coro}
\begin{example}
For $n=2$, consider $h:=E_{11}-E_{-1,-1}$, $e:=E_{12}-E_{-2,-1}$, and set
$h':=E_{22}-E_{-2,-2}$ as well as $f:=(1-et)^{-1}$.
By Corollary \ref{Jordanian}, we get a
Jordaninan quantization on
$\mathbf{u}_{t,q}(\mathfrak{sp}_{4})$ with the coproduct as
follows:
\begin{equation*}
\begin{split}
&\Delta(h)=1 \otimes h+h \otimes f,\\
&\Delta(h')=h' \otimes 1+1\otimes h' +h \otimes (1-f),\\
&\Delta(e)= 1\otimes e+e \otimes f^{-1},\\
\end{split}
\end{equation*}
\begin{equation*}
\begin{split}&\Delta(E_{1,-2}{+}E_{2,-1})=(E_{1,-2}{+}E_{2,-1}) \otimes f^{-1} +1
\otimes (E_{1,-2}{+}E_{2,-1})-2h
\otimes fE_{1,-1}t,\\
&\Delta(E_{-1,2}+E_{-2,1})=(E_{-1,2}+E_{-2,1}) \otimes f+1 \otimes
(E_{-1,2}+E_{-2,1})+2h \otimes fE_{-2,2}t,\\
&\Delta(-E_{-1,-2}+E_{2,1})=(-E_{-1,-2}+E_{2,1}) \otimes f+1 \otimes
(-E_{-1,-2}+E_{2,1})\\
&\hskip3cm-h \otimes f(h-h')t-h^{\langle 2\rangle}\otimes f^2et^2,\\
&\Delta(E_{1,-1})=E_{1,-1} \otimes f^{-2}+1 \otimes E_{1,-1},\\
&\Delta(E_{2,-2})=E_{2,-2} \otimes 1 + 1 \otimes E_{2,-2} - h
\otimes f(E_{1,-2}+E_{2,-1}) t+h^{\lg 2 \rg} \otimes
f^2 E_{1,-1} t^2,\\
&\Delta(E_{-1,1})=E_{-1,1} \otimes f^2 + 1 \otimes E_{-1,1}+h \otimes
f (E_{-1,2}+E_{-2,1})t+h^{\lg 2 \rg} \otimes f^2
E_{-2,2}t^2,\\
&\Delta(E_{-2,2})=E_{-2,2} \otimes 1 + 1 \otimes E_{-2,2}.
\end{split}
\end{equation*}
\end{example}

\medskip
\subsection{Open Questions} As is well known, Andruskiewitsch and Schneider \cite{AS}
gave a certain classification (i.e., a distinguished isomorphism theorem)
for the finite-dimensional complex pointed Hopf algebras with abelian finite group algebras
as their coradicals whose orders satisfy some conditions.
This is the most important achievement in Hopf algebra theory during the last two decades.
However, the new Hopf algebras of prime-power dimensions that we obtained above do not belong
to the class of pointed ones obtained over the complex numbers.

In conclusion, we would like to propose the following
interesting questions for further consideration.
\medskip

{\bf Question 1.} \ Assume that the ground field is algebraically closed of prime characteristic.
How many non-isomorphic Hopf algebra structures can be obtained
on the restricted universal enveloping algebras $\mathbf u(H(2n;\underline{1}))$ or
$\mathbf u(\mathfrak{sp}_{2n})$ by Drinfel'd twists? Are there infinitely many isomorphism classes as in the
famous counter-example to Kaplansky's 10th conjecture (see \cite{K}) obtained by
Beattie et al \cite{BDG}? Is it possible to classify the
Hopf algebras of a given prime-power dimension (up to isomorphism)?

\medskip
{\bf Question 2.} \ What are the conditions for $\mathbf
u_{t,q}(H(2n;\underline{1}))$ or $\mathbf u_{t,q}(\mathfrak {sp}_{2n})$ to be
a ribbon Hopf algebra (see \cite{HW2} and references therein)?

\medskip
{\bf Question 3.} \ It might be interesting to consider the
tensor product structures of representations for $\mathbf
u_{t,q}(H(2n;\underline{1}))$ or $\mathbf u_{t,q}(\mathfrak {sp}_{2n})$, respectively. How
do their tensor categories behave?

\bibliographystyle{amsalpha}

\end{document}